\documentclass{amsart}
\usepackage{righttag}
\usepackage{epsf}
\def\hepsffile{\leavevmode\epsffile}

\theoremstyle{plain}
\newtheorem{thm}{Theorem}[subsection]

\newtheorem{lem}[thm]{Lemma}

\theoremstyle{remark}
\newtheorem{rem}{Remark}
 
\theoremstyle{definition}
\newtheorem{defin}[thm]{Definition}

\newtheorem{emf}[thm]{}

\def\ind{\protect\operatorname{ind}}

\def\Int{\protect\operatorname{Int}}

\def\mod{\protect\operatorname{mod}}

\def\gl{\protect\operatorname{gl}}

\def\pr{\protect\operatorname{pr}}

\def\Ga{\alpha}

\def\C{{\mathbb C}}
\def\Z{{\mathbb Z}}

\def\R{{\mathbb R}}

\def\1{\hbox{\rm\rlap {1}\hskip.03in{\textrm I}}}
\def\Bbbone{{\rm1\mathchoice{\kern-0.25em}{\kern-0.25em}
	{\kern-0.2em}{\kern-0.2em}I}}

\def\p{\partial}

\def\hepsffile{\leavevmode\epsffile} 
\begin{document} 
\title[Invariants of fronts on surfaces]{Shadows of wave fronts\\
and\\
Arnold-Bennequin type invariants of fronts \\
on surfaces and orbifolds}
\author[V.~Tchernov]{Vladimir Tchernov}
\address{D-MATH, ETH Zentrum, CH-8092, Z\"urich, Switzerland}
\email{Chernov@math.ethz.ch}
%\date{\today}

\centerline{\em This article has been published 
in:} 
\centerline{\em Amer. Math. Soc. Transl. (2) Vol. {\bf 190}, 1999,
pp. 153-184.\/}

\begin{abstract}

A first order Vassiliev invariant of an oriented 
knot in an $S^1$-fibration and a Seifert fibration 
over a surface is constructed. 
It takes values in a quotient of the group ring of 
the first homology group of the total space of the fibration. 
It gives rise
to an invariant of wave fronts on surfaces and orbifolds
related to the Bennequin-type
invariants of the Legendrian curves studied by F.~Aicardi, 
V.~Arnold, M.~Polyak, and S.~Tabachnikov.
Formulas expressing these relations are presented.

We also calculate Turaev's shadow for the Legendrian lifting of
a wave front. This allows to use in the case of wave fronts 
all invariants known for shadows.

\end{abstract}
\maketitle

Most of the proofs in this paper are postponed until the last section.

Everywhere in this text an $S^1$-fibration is a locally-trivial fibration
with fibers homeomorphic to $S^1$. 

In this paper the multiplicative notation for the operation 
in the first homology group is used. The zero homology class is 
denoted by $e$. The reason for this is that we have to deal with the 
integer group ring of the
first homology group. For a group $G$ the group
of all formal half-integer linear combinations of elements of $G$ is denoted 
by $\frac{1}{2}\Z[G]$.

We work in the differential category. 

I am deeply grateful to Oleg Viro for the inspiration of this work and many
enlightening discussions. I am thankful to Francesca Aicardi, Thomas
Fiedler, and Michael Polyak for the valuable discussions.  

\section{Introduction}
In~\cite{Polyak} Polyak suggested a quantization $l_q(L)\in
\frac{1}{2}\Z[q,q^{-1}]$ of the Bennequin 
invariant of a
generic cooriented oriented wave front $L\subset\R^2$. 
In this paper we construct an 
invariant $S(L)$ which is,  in a sense, 
a generalization of $l_q(L)$ to the case of
a wave front on an arbitrary surface $F$.  

In the same paper~\cite{Polyak} Polyak 
introduced Arnold's~\cite{Arnold} $J^+$-type
invariant of a front $L$ on an oriented surface $F$. It takes values in
$H_1(ST^*F,\frac{1}{2}\Z)$. 
We show that $S(L)\in \frac{1}{2}\Z[H_1(ST^*F)]$ is a refinement of this
invariant in the sense that it is taken to 
Polyak's invariant under the natural mapping 
$\frac{1}{2}\Z[H_1(ST^*F)]\rightarrow H_1(ST^*F,\frac{1}{2}\Z)$. 

Further we generalize $S(L)$ to the case where $L$ is a wave front
on an orbifold.

Invariant $S(L)$ is constructed in two steps. The first one consists in 
lifting of $L$ to the Legendrian knot
$\lambda$ in the $S^1$-fibration
$\pi:ST^*F\rightarrow F$. The second step can be applied to any knot in an
$S^1$-fibration, and 
it involves the structure of the fibration in a
crucial way. This step allows us to define the
$S_K$ invariant of a knot $K$ in the total space $N$ of an 
$S^1$-fibration. Since ordinary knots are considered up to a rougher
equivalence relation (ordinary isotopy versus Legendrian isotopy), in order
for $S_K$ to be well defined it has to take values in a quotient of 
$\Z[H_1(N)]$. This invariant is
generalized to the case of a knot in a Seifert fibration, and this allows
us to
define $S(L)$ for wave fronts on orbifolds.

All these invariants are Vassiliev invariants of order one in an
appropriate sense. 

For each of these invariants we introduce its version with values in 
the group of formal linear combinations of the 
free homotopy classes of oriented curves in the
total space of the corresponding fibration.

The first invariants of this kind were constructed by
Fiedler~\cite{Fiedler} in the case
of a knot $K$ in a $\R^1$-fibration over a surface and by Aicardi in the
case of
a generic oriented cooriented wave front $L\subset \R^2$. 
The connection between
these invariants and $S_K$ is discussed in~\cite{Tchernov}.

The space $ST^*F$ is naturally fibered over a surface $F$  with a fiber
$S^1$. In~\cite{Turaev} Turaev introduced a shadow description 
of a knot $K$ in
an oriented three dimensional manifold $N$ fibered over an oriented surface 
with a fiber 
$S^1$. A shadow presentation of a knot $K$ is a generic projection of $K$,
together with an assignment of numbers to regions. It describes a knot type
modulo a natural action of $H_1(F)$.
It appeared to be a very useful
tool. Many invariants of knots in $S^1$-fibrations, in particular quantum
state sums, can be expressed
as state sums for their shadows. In this paper
we construct shadows of Legendrian liftings of wave fronts. 
This allows one to use any
invariant already known for shadows in the case of wave fronts. 

However, in this paper shadows are used mainly for the purpose of depicting
knots in $S^1$-fibrations.

\section{Shadows}\label{sh-def}
\subsection{Preliminary constructions}\label{prelim}

We say that a one-dimensional submanifold $L$ of a total space $N^3$ of 
a fibration 
$\pi:N^3\rightarrow F^2$ is {\em generic
with respect to $\pi$\/ } if $\pi\big|_L$ is a generic immersion.
Recall that an immersion of 1-manifold into a surface is said to
be {\em generic\/} if it has neither self-intersection points with
multiplicity greater than $2$ nor self-tangency  points, and at
each double point its branches are transversal to each other. An
immersion of (a circle) $S^1$ to a surface is called a {\em curve\/}.

Let $\pi$ be an oriented $S^1$-fibration of $N$ over an oriented
closed surface  $F$.

$N$ admits a fixed point free involution which preserves fibers.
Let $\tilde N$ be the quotient of $N$ by
this involution, and let $p:N\rightarrow\tilde{N}$ be the corresponding double
covering. Each fiber of $p$ (a pair of antipodal points) is
contained in a fiber of $\pi$. Therefore, $\pi$ factorizes through
$p$ and we have a fibration $\tilde{\pi}:\tilde{N}\rightarrow F$.
Fibers of $\tilde\pi$ are projective lines. They are
homeomorphic to circles.

An isotopy of a link $L\subset N$ is said to be vertical with respect to
$\pi$ if each point of $L$ moves along a fiber of $\pi$. It is clear
that if two links are vertically isotopic, then their projections
coincide. Using vertical isotopy we can modify each generic
link $L$ in such a way that any two points of $L$ belonging to the
same fiber lie in the same orbit of the involution. 
Denote the obtained generic link by $L'$.

Let $\tilde{L}=p(L')$. It is obtained from $L'$ by gluing together
points lying over the same point of $F$. Hence $\tilde{\pi}$ maps
$\tilde{L}$ bijectively to $\pi(L)=\pi(L')$. Let $r:\pi (L)\rightarrow
\tilde{L}$ be an inverse bijection. It is a section of
$\tilde{\pi}$ over $\pi (L)$.

For a generic non-empty collection of curves on a surface by a {\em region\/}
we mean the closure of a connected component of the complement of
this collection. Let $X$ be a region for $\pi (L)$ on $F$. Then 
$\tilde{\pi}\big|_X$ is a trivial fibration. Hence we can identify it with
the projection $S^1 \times X\rightarrow X$. Let $\phi$ be a
composition of the section $r\big|_{\p X}$ with the projection to $S^1$.
It maps $\p X$ to  $S^1$.
Denote by $\Ga_X$ the degree of $\phi$. (This is actually an
obstruction to an extension of $r\big|_{\p X}$ to $X$.) One can
see that $\Ga_X$ does not depend on the choice of the trivialization
of $\tilde{\pi}$ and on the choice of $L'$.

\subsection{Basic definitions and properties}
\begin{defin}\label{def-shadowlink}
The number $\frac{1}{2}\alpha _X$ corresponding to a region $X$ is
called the {\em gleam\/} of $X$ and is denoted by $\gl(X)$.
A {\em shadow\ $s(L)$ of a generic link $L\subset N$\/} is a (generic)
collection of curves $\pi (L)\subset F$ with the gleams assigned to each
region $X$.  The sum of
gleams over all regions is said to be the {\em total gleam\/} of the
shadow.  
\end{defin}

\begin{emf}\label{prop1-gleam}
One can check that for any region $X$ the integer $\alpha _X$ is congruent
modulo 2 to the number of corners of $X$. Therefore, $\gl(X)$ is an integer if
the region $X$ has even number of corners and half-integer otherwise.
\end{emf}

\begin{emf}\label{prop2-gleam}
The total gleam of the shadow is equal to the Euler number of $\pi$.
\end{emf}

\begin{defin}\label{def-shadow}
A {\em shadow \/} on $F$ is a generic collection of curves together
with the numbers $\gl(X)$ assigned to each region $X$. These numbers
can be either integers or half-integers, and they should satisfy the
conditions of~\ref{prop1-gleam}~and~\ref{prop2-gleam}.
\end{defin}

\begin{figure}[htbp]
 \begin{center}
  \epsfxsize 10cm
  \hepsffile{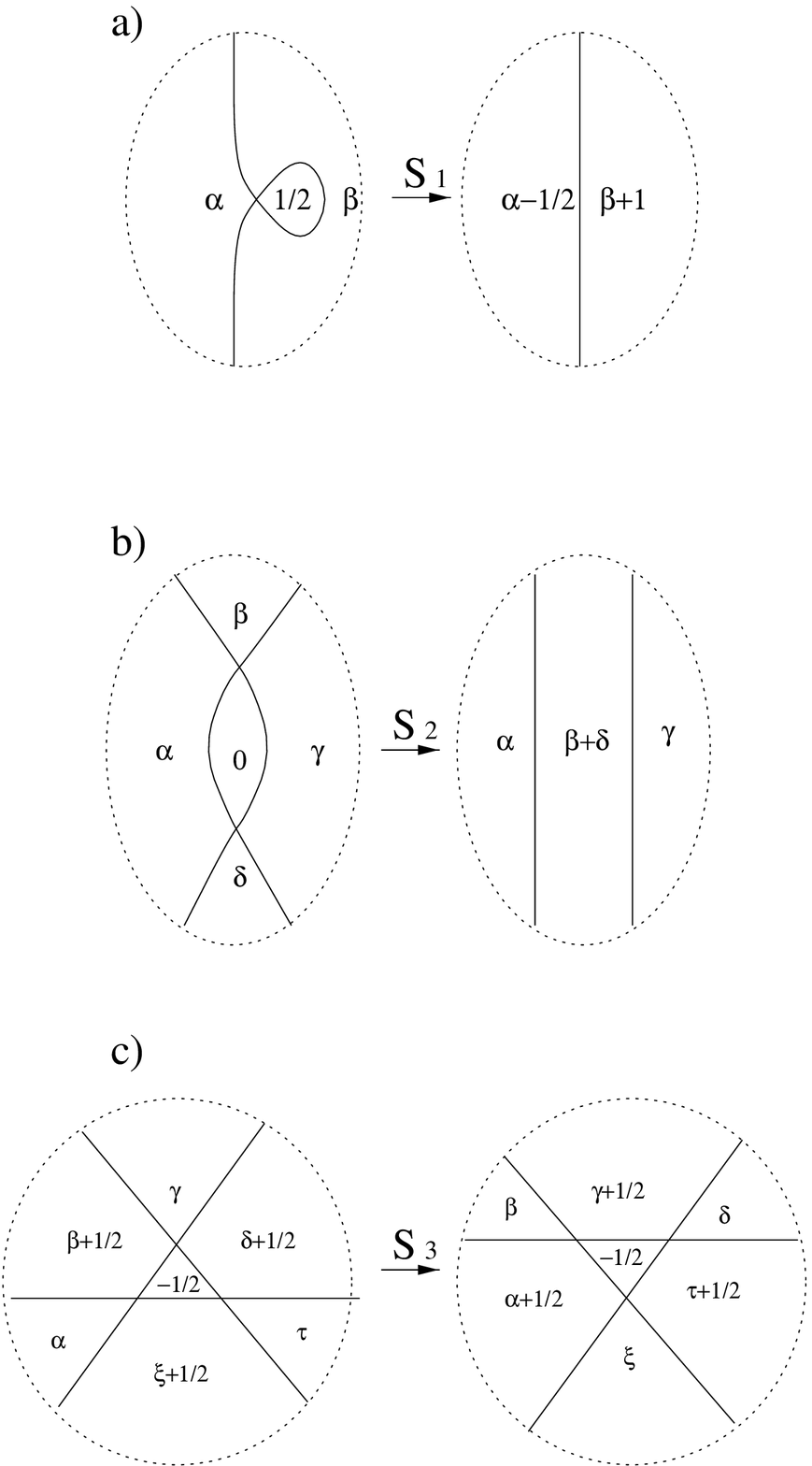}
 \end{center}
\caption{}\label{shad3.fig}
\end{figure}

There are three local moves $S_1,S_2$, and $S_3$ of shadows shown in
Figure~\ref{shad3.fig}. They are similar to the well-known
Riedemeister moves of planar knot diagrams.

\begin{defin} Two shadows are said to be {\em shadow equivalent\/} if
they can be transformed to each other by a finite sequence of moves
$S_1,S_2,S_3$, and their inverses.
\end{defin}

\begin{figure}[htbp]
 \begin{center}
  \epsfxsize 10cm
  \hepsffile{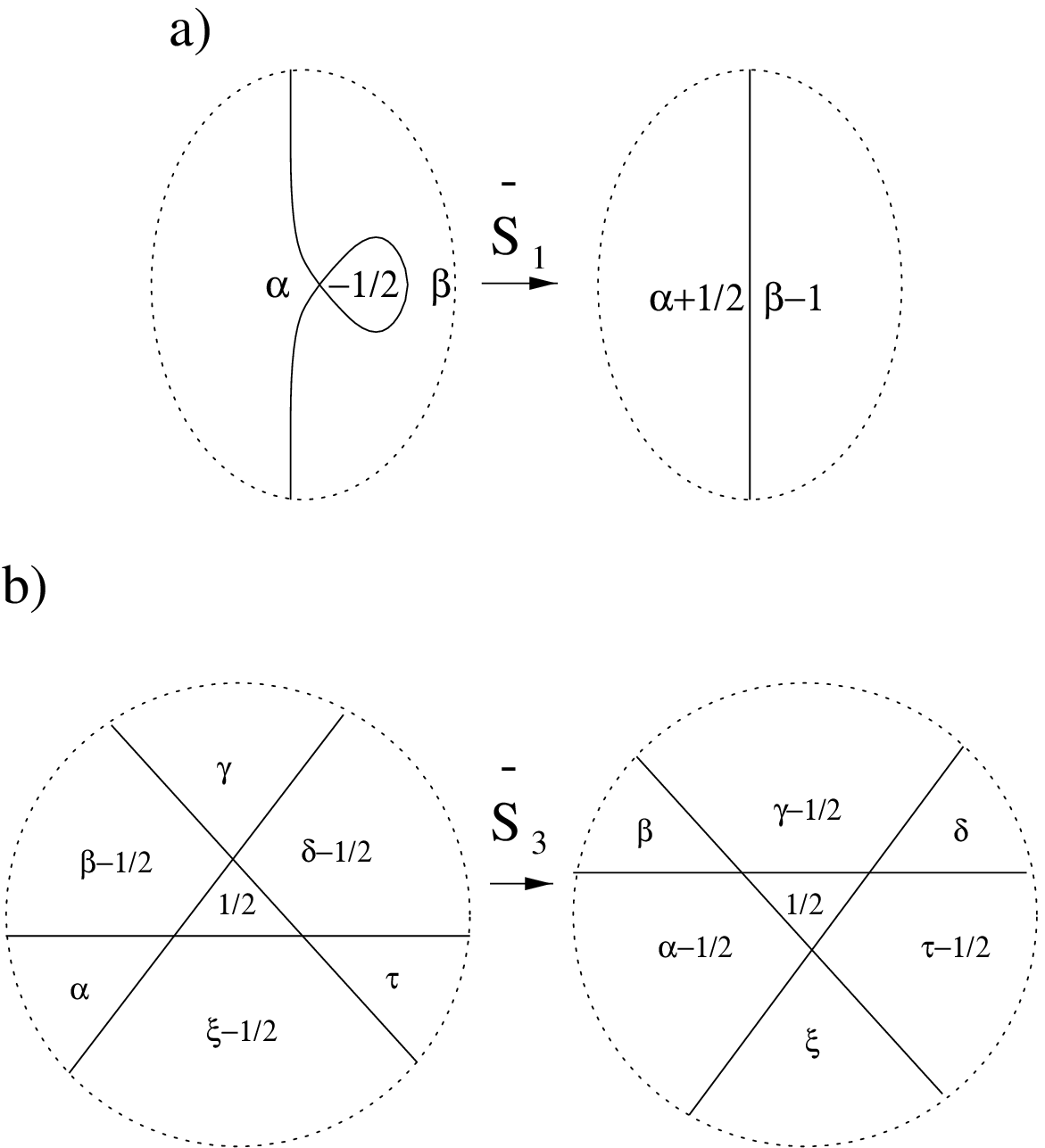}
 \end{center}
\caption{}\label{shad5.fig}
\end{figure}

\begin{emf}
 There are two more important shadow moves $\bar S_1$ and $\bar S_3$ 
shown in Figure~\ref{shad5.fig}. They are similar to the previous versions of 
the first and the third Riedemeister
moves and can be expressed 
in terms of $S_1,S_2$, $S_3$, and their inverses.
\end{emf}

\begin{emf}\label{actionconstruct}
In~\cite{Turaev} the action of 
$H_1(F)$ on the set of all isotopy
types of links in $N$ is constructed as follows. 
Let $L$ be a generic link in $N$ and $\beta$ an
oriented (possibly self-intersecting) curve on $F$ presenting
a homology class $[\beta ]\in H_1(F)$. Deforming $\beta $ we can
assume that $\beta $ intersects $\pi(K)$ transversally at a finite
number of points different from the self-intersection points of $\pi(K)$.
Denote by $\alpha =[a,b]$ a small segment of $L$ such that
$\pi(\alpha )$  contains exactly one intersection point $c$ of $\pi(L)$
and $\beta $. Assume that $\pi(a)$ lies to the left, and
$\pi(b)$ to the right of $\beta $. Replace $\alpha $ by the arc
$\alpha '$ shown in Figure~\ref{shad13.fig}.
\begin{figure}[htbp]
\begin{center}
\epsfxsize 5cm
  \hepsffile{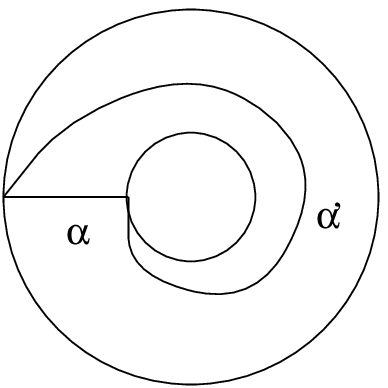}
 \end{center}
\caption{}\label{shad13.fig}
\end{figure}
 We will call this
transformation of $L$ a {\em fiber fusion\/} over the point $c$.
After we apply fiber fusion to $L$ over all points of $\pi(L)\cap\beta$
we get a new generic link $L'$ with $\pi(L)=\pi(L')$. One can notice that
the shadows of $K$ and $K'$ coincide. Indeed, each time  $\beta $
enters a region $X$ of $s(L)$, it must leave it. 
Hence the contributions of the newly inserted arcs to
the gleam of $X$ cancel out. Thus links belonging to one
$H_1(F)$-orbit always produce the same shadow-link on $F$.

\end{emf}

\begin{thm}[Turaev \cite{Turaev}]\label{action}
Let $N$ be an oriented closed manifold, $F$ an oriented surface, and 
$\pi:N\rightarrow F$ an
$S^1$-fibration with the Euler number $\chi (\pi)$. The
mapping that associates to each link $L\subset N$ its shadow equivalence
class on $F$ establishes a bijective correspondence between the
set of isotopy types of links in $N$ modulo the action of
$H_1(F)$ and the set of all shadow equivalence classes on $F$ 
with the total gleam $\chi(\pi)$. 
\end{thm}

\begin{emf}\label{homolfusion}
It is easy to see that all links whose projections represent 
$0\in H_1(F)$ and whose shadows coincide are homologous 
to each other. To prove this,
one looks at the description of a fiber fusion and notices that to each
fiber fusion where we add a positive fiber corresponds another  where we
add a negative one. Thus the numbers of positively and negatively 
oriented fibers we add are equal, and they cancel
out.
\end{emf}

\begin{emf}\label{shadowgeneral}
As it was remarked in~\cite{Turaev} it is easy to transfer the construction 
of shadows and Theorem~\ref{action} to the 
case where $F$ is a non-closed oriented surface and $N$ is an oriented
manifold. In order
to define the gleams of the regions that have a non-compact closure or
contain components of $\p F$, we have to choose a section of the fibration 
over all boundary components and ends of $F$. In the case of non-closed $F$
the total gleam of the shadow is equal to the obstruction 
to the extension of the section to the entire surface.
\end{emf}

\section{Invariants of knots in $S^1$-fibrations.}
\subsection{Main constructions}

In this section we deal with knots in an $S^1$-fibration
$\pi$ of an oriented three-dimensional manifold $N$ over an oriented 
surface $F$. We do not assume $F$ and $N$ to be
closed. As it was said in~\ref{shadowgeneral}, all theorems 
from the previous section are applicable in this case.

\begin{defin}[of $S_K$]\label{SK} Orientations of $N$ and $F$ determine an 
orientation of a fiber of the fibration. Denote by $f\in H_1(N)$ 
the homology class of a positively oriented fiber.

Let $K\subset N$ be an oriented knot which is generic with respect to $\pi$.
Let $v$ be a double point of $\pi (K)$. The fiber $\pi^{-1}(v)$ divides $K$ 
into two arcs that inherit the orientation from $K$. Complete each arc
of $K$ to an oriented knot by adding the arc of $\pi^{-1}(v)$ 
such that the orientations of these
two arcs define an orientation of their union. The orientations of $F$ and
$\pi(K)$ allow one to identify a small neighborhood of $v$ in $F$ with a
model picture shown in Figure~\ref{shad1.fig}a. Denote the knots
obtained by the operation above by $\mu^+_v$ and $\mu^-_v$ as shown in
Figure~\ref{shad1.fig}. We will often call this construction a {\em
splitting\/} of $K$ (with respect to the orientation of $K$).

This splitting can be described in terms of shadows as follows.
Note that $\mu ^+_v$ and $\mu ^+_v$ are not in general position. We slightly 
deform them in a neighborhood of $\pi^{-1}(v)$, so that
$\pi(\mu ^+_v)$ and $\pi(\mu ^+_v)$ 
do not have double points in the neighborhood of $v$.
Let $P$ be a neighborhood of $v$ in $F$
homeomorphic to a closed disk.  
Fix a section over $\p P$ such that the intersection
points of $K\cap\pi^{-1}(\p P)$ belong to the section. Inside $P$ we can
construct Turaev's shadow (see \ref{shadowgeneral}). 
The action  of $H_1(\Int P)=e$  
on the set of the isotopy types of links is trivial. 
Thus the part of $K$ can be reconstructed in the unique way 
(up to an isotopy fixed on $\p P$) from the shadow over $P$
(see~\ref{shadowgeneral}). 
The shadows for $\mu^+_v$ and $\mu^-_v$ are 
shown in Figures~\ref{shad1.fig}a~and~\ref{shad1.fig}b respectively.

\begin{figure}[htb]
\begin{center}
 \epsfxsize 10 cm
\hepsffile{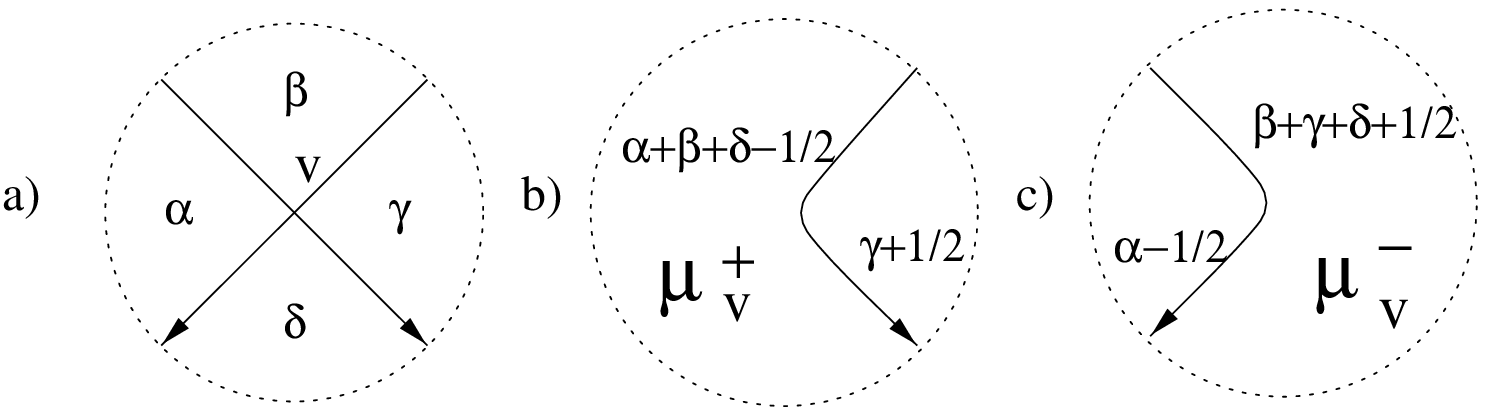}
 \end{center}
\caption{}\label{shad1.fig}
\end{figure}

Regions for the shadows $s(\mu ^+_v)$ and 
$s(\mu ^-_v)$ are, in fact, unions of regions for $s(K)$. One should think of 
gleams as of measure, so that the gleam of a region is the sum 
of all numbers inside.

Let $H$ be the quotient of the group ring $\Z [H_1(N)]$ (viewed as a
$\Z$-module) 
by the submodule generated by $\bigl\{ [K]-f, [K]f-e \bigr\} $. 
Here by $[K]\in H_1(N)$ we denote the homology class represented by the image
of $K$.

Finally define $S_K\in H$ by the following
formula, where the summation is taken over all double points $v$
of $\pi (K)$:

\begin{equation}\label{eqSK}
S_K=\sum_v \bigl([\mu ^+_v]-[\mu ^-_v]\bigr).
\end{equation}
\end{defin}

\begin{emf}\label{sum-property} Since $\mu^+_v\cup\mu^-_v=K\cup\pi^{-1}(v)$
we have
\begin{equation}
[\mu ^+_v][\mu ^-_v]= [K]f.
\end{equation}
\end{emf}

\begin{thm}\label{correct2} $S_K$ is an isotopy invariant of the knot $K$.
\end{thm}

For the proof of Theorem~\ref{correct2} see Subsection~\ref{pfcorrect2}.

\begin{emf}\label{stupid}
It follows from~\ref{sum-property}  
that $S_K$ can also be described as an element of $\Z[H_1(N)]$ equal to the
sum of $\bigl([\mu^+_v]-[\mu^-_v]\bigr)$ over all double points for 
which the sets $\{ [\mu ^+_v],[\mu ^-_v] \}$ and $\{ e,f \}$ are disjoint. 
Note
that in this case we do not need to factorize $\Z[H_1(N)]$ to make $S_K$
well defined.
\end{emf}

\begin{emf}\label{homotopvers}
One can obtain an invariant similar to $S_K$ with
values in the free $\Z$-module generated by the set of all free 
homotopy classes of oriented curves in $N$. To do this one substitutes
the homology classes of $\mu ^+_v$ and $\mu ^-_v$ in~\eqref{eqSK} 
with their free homotopy
classes and takes the summation over the set of all double points $v$ of 
$\pi (K)$ such that neither one of the knots $\mu ^+_v$ and $\mu ^-_v$ 
is homotopic to
a trivial loop and neither one of them is homotopic to a positively
oriented fiber
(see~\ref{stupid}).

To prove that this is indeed an invariant of $K$ one can easily modify the
proof of Theorem~\ref{correct2}.   
\end{emf}

\subsection{$S_K$ is a Vassiliev invariant of order one}

\begin{emf} Let $\pi:N\rightarrow F$ be an $S^1$-fibration over a
surface. Let $K\subset N$ be a knot generic with respect to $\pi$ and $v$ a
double point of $\pi(K)$. A modification of pushing of one branch of $K$ 
through the other along the fiber $\pi ^{-1}(v)$ 
is called a {\em modification of $K$ along the fiber\/} $\pi^{-1}(v)$.
\end{emf}

\begin{emf}\label{helpdelta}
If a fiber fusion increases by one the gleam $\gamma$ 
in Figure~\ref{shad1.fig}b, then $[\mu ^+_v]$ is multiplied by $f$. 
If a fiber fusion increases by one the gleam $\alpha$ 
in Figure~\ref{shad1.fig}c, then $[\mu ^-_v]$ is multiplied by $f^{-1}$. 
These facts are easy to verify.
\end{emf}

\begin{emf}\label{vassiliev2}   
Let us find out how $S_K$ changes under the modification along a fiber over a
double point $v$. Consider a singular knot $K'$
(whose only singularity is a point $v$ of transverse
self-intersection). Let $\xi_1$ and $\xi_2$ be the homology classes of the
two loops of $K'$ adjacent to $v$. The two
resolutions of this double point correspond to adding $\pm\frac{1}{2}$ to
the gleams of the regions adjacent to $v$ in two ways shown in
Figures~\ref{shad14.fig}b and~\ref{shad14.fig}c. 
\begin{figure}[htbp]
 \begin{center}
  \epsfxsize 12cm
  \hepsffile{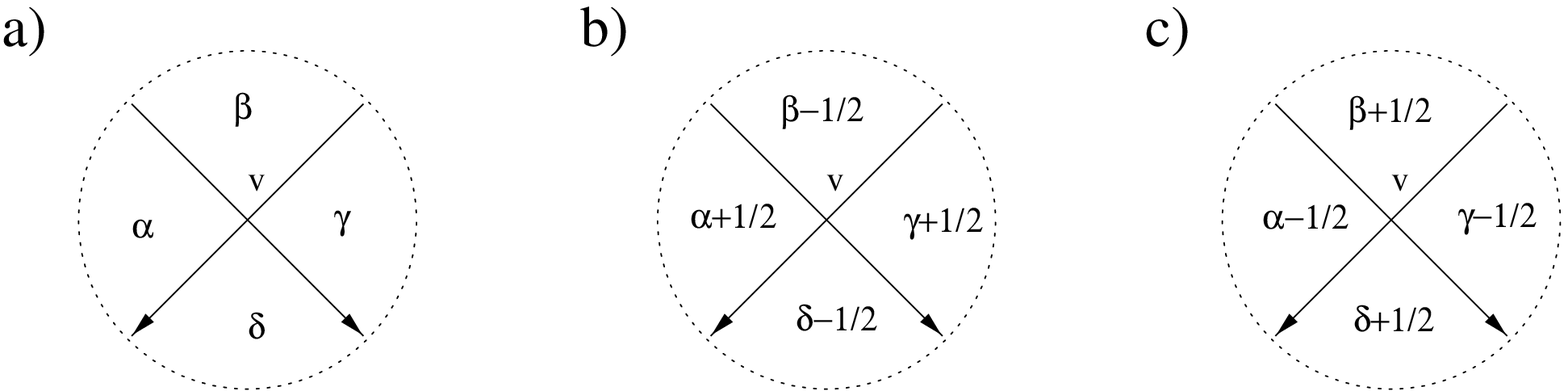}
 \end{center}
\caption{}\label{shad14.fig}
\end{figure} 

Using~\ref{helpdelta} one verifies that under the corresponding modification 
$S_K$ changes by 
\begin{equation}\label{type2}
(f-e)(\xi_1+\xi_2).
\end{equation}

This means that the first derivative of $S_K$ depends
only on the homology classes of the two loops adjacent to the singular
point. Hence the second derivative of $S_K$ is $0$. Thus it is a Vassiliev 
invariant of order one in the usual sense.

For similar reasons
the version of $S_K$ with values in the free $\Z$-module generated 
by all free homotopy classes of oriented curves in $N$ is also a Vassiliev
invariant of order one.
\end{emf}

\begin{thm}\label{realization2}$ $
\begin{description}
\item[\textrm{I}] If $K$ and $K'$ are two knots representing the same free
homotopy class, then $S_{K}$ and $S_{K'}$ are congruent modulo 
the submodule generated by elements of form
\begin{equation}\label{type}
(f-e)(j+[K]j^{-1})
\end{equation}
for $j\in H_1(N)$.

\item[\textrm{II}] If $K$ is a knot, and $S\in H$ is congruent to $S_{K}$ modulo the 
submodule generated by elements of 
form~\eqref{type} (for $j\in H_1(N)$), 
then there exists a knot $K'$ such that:

\begin{description}
\item[a] $K$ and $K'$ represent the same free homotopy class;

\item[b] $S_{K'}=S$.
\end{description}
\end{description}
\end{thm}

For the proof of Theorem~\ref{realization2} see
subsection~\ref{pfrealization2}.

\subsection{Example.}
If $N$ is a solid torus $T$ fibered over a disk, 
then we can calculate the value
of $S_K$ directly from the shadow of $K$.

\begin{figure}[htb]
\begin{center}
 \epsfxsize\hsize\advance\epsfxsize -0.5cm
 \hepsffile{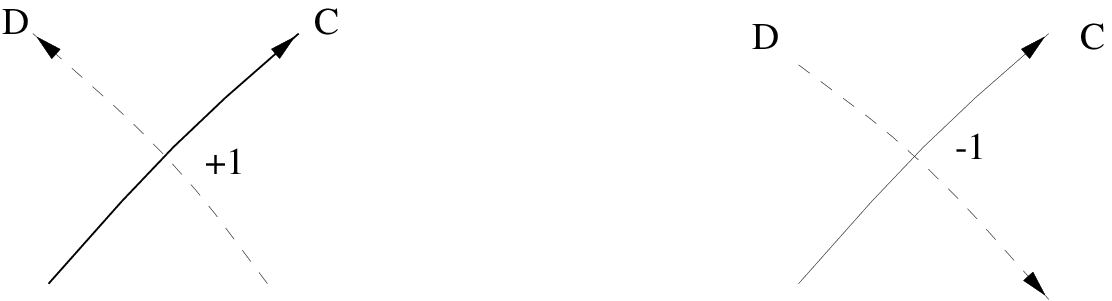}
 \end{center}
\caption{}\label{shad12.fig}
\end{figure}

\begin{defin} Let $C$ be an oriented closed curve in $\R^2$ and $X$ a
region for $C$. Take a point $x\in\Int X$ and connect it to a point near
infinity by a generic oriented path $D$. 
Define the sign of an intersection point of $C$ and $D$ as 
shown in Figure~\ref{shad12.fig}. Let $\ind_C X$ be the sum over
all intersection points of $C$ and $D$ of the signs of these points.
\end{defin}

It is easy to see that $\ind_C(X)$ is independent on the choices 
of $x$ and $D$.

\begin{defin} Let $K\subset T$ be an oriented knot which is
generic with respect to $\pi$, and let $s(K)$ be its
shadow. Define $\sigma(s(K))\in \Z$ as the following sum 
over all regions $X$ for $\pi (K)$: 
\begin{equation}
\sigma(s(K))=\sum_X\ind_{\pi(K)}(X)\gl(X).
\end{equation}
\end{defin}

Denote by $h\in \Z$ the image of $[K]$ under the natural identification of 
$H_1(T)$ with $\Z$. 

\begin{lem}\label{homtor}
$\sigma(s(K))=h.$ 
\end{lem}

\begin{emf}Put 

\begin{equation}
S'_K=\sum t^{\sigma(s(\mu ^+_v))}-t^{\sigma(s(\mu ^-_v))},
\end{equation}
where the sum is taken over all double points $v$ of $\pi (K)$
such that $\{ 0,1\}$ 
and $\{ \sigma(s(\mu ^+_v)),\sigma(s(\mu ^-_v))\}$  are disjoint
(see~\ref{stupid}).
   
Lemma~\ref{homtor} implies that $S'_K$ is the image of $S_K$ under the
natural identification of $\Z[H_1(T)]$ with the ring of
finite Laurent polynomials (see~\ref{stupid}).

\end{emf}

One can show~\cite{Tchernov} that $S'_K$ and Aicardi's partial
linking polynomial of $K$ (which was introduced in~\cite{Aicardi})
can be explicitly expressed in terms of each other.

\subsection{Further generalizations of the $S_K$ invariant}
One can show that an invariant similar to $S_K$ can be introduced in the
case where $N$ is oriented and $F$ is non-orientable.

\begin{figure}[htb]
 \begin{center}
  \epsfxsize 12cm
  \hepsffile{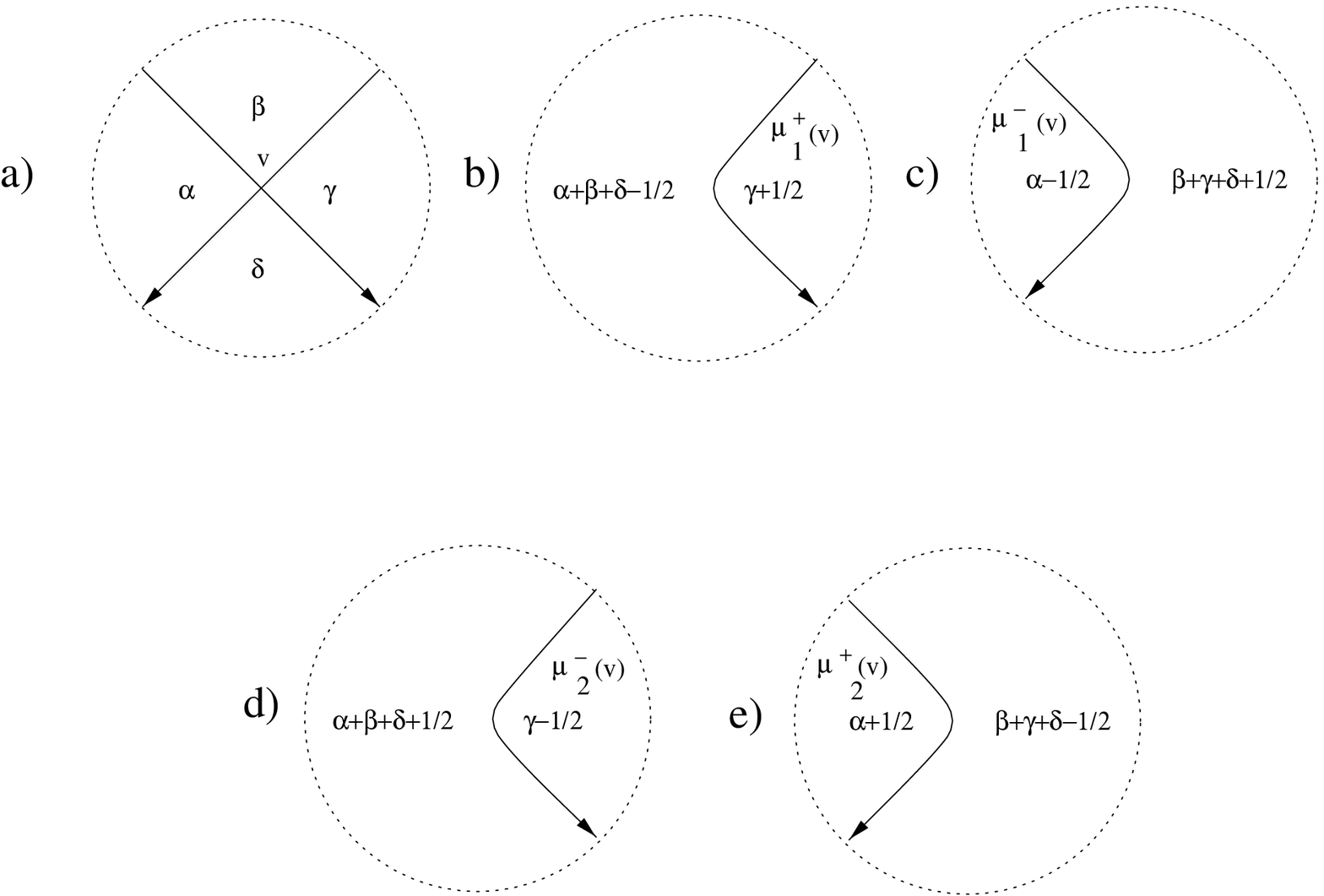}
 \end{center}
\caption{}\label{shadun2.fig}
\end{figure}

\begin{defin}[of $\tilde S_K$]\label{tildeSK}
Let $N$ be oriented and $F$ non-orientable.
Let $K\subset N$ be an oriented knot generic with respect to $\pi$, and 
let $v$ be a double point of $\pi(K)$. Fix an orientation of a small
neighborhood of $v$ in $F$. Since $N$ is oriented
this induces an orientation
of the fiber $\pi^{-1}(v)$. Similarly to the definition of $S_K$
(see~\ref{SK}), we split our knot with respect to the 
orientation and obtain two
knots $\mu_1^+(v)$ and $\mu_1^-(v)$. Then we take the other orientation 
of the neighborhood of $v$ in $F$, and in the same way we obtain another pair
of knots $\mu_2^+(v)$ and $\mu_2^-(v)$. The element 
$\bigl( [\mu_1^+(v)]-[\mu_1^-(v)]+[\mu_2^+(v)]-[\mu_2^-(v)]\bigr)\in
\Z[H_1(N)]$ does not depend on which orientation of the neighborhood of $v$
we choose first. 

Similarly to the definition of $S_K$, we can describe all
this in terms of shadows as it is shown in Figure~\ref{shadun2.fig}.
These shadows are constructed with respect to the same orientation of the
neighborhood of $v$.

Let $f$ be the homology class of a fiber of $\pi$ oriented in some way. 
As one can easily prove $f^2=e$, so it
does not matter which orientation we choose to define $f$.
Let $\tilde H$ be the quotient of $\Z [H_1(N)]$ (viewed as a $\Z$-module) by the 
$\Z$-submodule
generated by $\Bigl \{ [K]-f+e-[K]f=(e-f)([K]+e) \Bigr\}$.
Finally define $\tilde S_K\in \tilde H$ by the
following formula, where the summation is taken 
over all double points
$v$ of $\pi (K)$: 
\begin{equation}\label{eqtildeSK}
\tilde 
S_K=\sum_v \Bigl([\mu _1^+(v)]-[\mu _1^-(v)]+[\mu _2^+(v)]-[\mu
_2^-(v)]\Bigr).
\end{equation}
\end{defin}

\begin{thm}\label{correct3} $\tilde 
S_K$ is an isotopy invariant of the knot $K$.
\end{thm}

The proof is essentially the same as the proof of Theorem~\ref{correct2}.

\begin{emf} One can easily prove that $\tilde S_K$ 
invariant satisfies
relations similar to~\eqref{type2}. In particular, $\tilde S_K$ is also a
Vassiliev invariant of order one. 
%The analogue of Theorem~\ref{realization2}
%also holds for $\tilde S_K$. 

One can introduce a version of this invariant with values in the free
$\Z$-module generated by all free homotopy classes of oriented curves in 
$N$. 
To do this, we substitute the homology classes of $\mu _1^+(v),$
$\mu _1^-(v),$ $\mu _2^+(v)$, and $\mu_2^-(v)$ with the corresponding 
free homotopy classes. The summation should be taken over the set of all
double points of $\pi(K)$ for which neither one of $\mu _1^+(v),
\mu _1^-(v),\mu _2^+(v)$, and $\mu_2^-(v)$ is homotopic to a trivial loop and
neither one of them is homotopic to a fiber of $\pi$. To prove that this is
indeed an invariant of $K$, one easily modifies the proof of
Theorem~\ref{correct2}. 
\end{emf}

\section{Invariants of knots in Seifert fibered spaces}

Let $(\mu, \nu)$ be a pair of relatively prime integers. Let 
$$D^2=\Bigl\{ (r,\theta); 0\leq r\leq 1, 0\leq \theta \leq 2\pi\Bigr\}\subset
\R ^2$$
be the unit disk defined in polar coordinates. A fibered solid torus of 
type $(\mu, \nu )$ is the quotient space of the cylinder $D^2\times I$ via
the identification $\bigl (\bigl (r,\theta \bigr ),1\bigr )=
\bigl (\bigl (r,\theta+\frac{2\pi \nu}{\mu}\bigr
),0\bigr)$. The fibers are the images of the curves ${x}\times I$.  The
integer $\mu$ is called the index or the multiplicity. For $|\mu| >1$ the
fibered solid torus is said to be {\em exceptionally fibered,\/} and the
fiber that is the image of $0\times I$ is called the {\em exceptional fiber\/}. 
Otherwise the
fibered solid torus is said to be {\em regularly fibered,\/} 
and each fiber is a {\em regular fiber\/}.

\begin{defin}\label{Seifert}
An orientable three manifold $S$ is said to be a 
{\em Seifert fibered manifold\/}
if it is a union of pairwise disjoint closed curves, called fibers, such
that each one has a closed neighborhood which is a union of fibers and
is homeomorphic to a fibered solid torus by a fiber preserving
homeomorphism.
\end{defin}

A fiber $h$ is called {\em exceptional\/} if $h$ has a neighborhood 
homeomorphic to
an exceptionally fibered solid torus (by a fiber preserving homeomorphism), 
and $h$ corresponds via the
homeomorphism to the exceptional fiber of the solid torus. If $\p S\neq
\emptyset$,
then $\p S$ should be a union of regular fibers. 

The quotient space obtained from a Seifert fibered manifold $S$ by
identifying each fiber to a point is a 2-manifold. It is called the orbit
space and the images of the exceptional fibers are called
{\em the cone points.\/}

\begin{emf}\label{numbers}
For an exceptional fiber $a$ of
an oriented Seifert fibered manifold there is a unique pair of relatively
prime integers $(\mu_a, \nu_a)$ such that $\mu_a>0$, $|\nu_a|<\mu_a$, and a
neighborhood of $a$ 
is homeomorphic (by a fiber preserving homeomorphism) to a
fibered solid torus of type $(\mu_a, \nu_a)$. We call the pair $(\mu_a,
\nu_a)$ {\em the type of the exceptional fiber\/} $a$. We also call this
pair the type of the corresponding cone point. 
\end{emf}

We can define an invariant of an oriented knot
in a Seifert fibered manifold that is similar to the $S_K$ invariant.

Clearly any $S^{1}$-fibration can be viewed as a Seifert fibration without
cone points. This justifies the notation in the definition below.

\begin{defin}[of $S_K$]\label{SeifertSK}
Let $N$ be an oriented Seifert fibered manifold with an oriented orbit space 
$F$. Let $\pi :N\rightarrow F$ be the corresponding fibration and $K\subset
N$  an oriented knot in general position with respect to $\pi$. Assume
also that $K$ does not intersect the exceptional fibers. For each double
point $v$ of $\pi (K)$ we split $K$ into $\mu ^+_v$ and $\mu ^-_v$ (see
~\ref{SK}). Let $A$ be the set of all exceptional
fibers. Since $N$ and $F$ are oriented, we have an induced orientation of
each exceptional fiber $a\in A$. For $a\in A$ set $f_a$ to be the homology 
class of the fiber with this orientation. For $a\in A$ of type $(\mu_a,
\nu_a)$ (see~\ref{numbers})
set $N_1(a)=\bigl\{ k\in \{1,\dots,\mu_a\}| \frac{2\pi k\nu_a}{\mu_a}\text{ } mod
\text{ }2\pi\in (0,\pi]\bigr\}$, 
$N_2(a)=\bigl\{ k\in \{1,\dots,\mu_a\}| \frac{2\pi k\nu_a}{\mu_a}\text{ } mod
\text{ }2\pi\in (0,\pi)\bigr\}$. Define $R^1_a, R^2_a\in \Z[H_1(N)]$ by the
following formulas:

\begin{equation}
R^1_a=\sum_{k\in N_1(a)}\bigl([K]f_a^{\mu_a-k}-f_a^k\bigr)-
\sum_{k\in N_2(a)}\bigl(f_a^{\mu_a-k}-[K]f_a^k\bigr),
\end{equation}
\begin{equation}
R^2_a=\sum_{k\in N_1(a)}\bigl(f_a^{\mu_a-k}-[K]f_a^k\bigr)-
\sum_{k\in N_2(a)}\bigl([K]f_a^{\mu_a-k}-f_a^k\bigr).
\end{equation}

Let $H$ be the quotient of $\Z [(H_1(N)]$ (viewed as a $\Z$-module) 
by the free $\Z$-submodule 
generated by 
$\Bigl \{[K]f-e, [K]-f, \bigl \{R^1_a, R^2_a\bigr\} _{a\in A}\Bigr \}$.
Finally, define $S_K\in H$ by the following formula, where the 
summation is taken over all double points $v$ of
$\pi (K)$:

\begin{equation}\label{SSK}
S_K=\sum_v\bigl( [\mu ^+_v]-[\mu ^-_v]\bigr).
\end{equation}
  
\end{defin}

\begin{thm}\label{correctSSK}
$S_K$ is an isotopy invariant of the knot $K$.
\end{thm}

For the proof of Theorem~\ref{correctSSK} see Subsection~\ref{correctSSK}.

We introduce a similar invariant in the case 
where $N$ is oriented and $F$ is non-orientable.

\begin{defin}[of $\tilde S_K$]\label{tilde1SK}

Let $N$ be an oriented Seifert fibered manifold with a non-orientable orbit
space $F$. Let $\pi:N\rightarrow F$ be the corresponding fibration and
$K\subset N$ an oriented knot in general position with respect to $\pi$.
Assume also that $K$ does not intersect the exceptional fibers. For each
double point $v$ of $\pi(K)$ we split $K$ into 
$\mu_1^+(v),$ $\mu_1^-(v),$ $\mu_2^+(v)$, and $\mu_2^-(v)$ 
as in~\ref{tildeSK}. The element 
$\bigl([\mu_1^+(v)]-[\mu_1^-(v)]+[\mu_2^+(v)]-[\mu_2^-(v)]\bigr)\in
\Z[H_1(N)]$ 
is well
defined. 

Denote by $f$ the homology class of a regular fiber oriented in some way. 
Note that $f^2=e$, so the orientation we use to define $f$ does not matter. 
For a cone point $a$ denote by $f_a$ 
the homology class of the fiber $\pi^{-1}(a)$ oriented in some
way.
   
For $a\in A$ of type $(\mu_a,\nu_a)$ set 
$N_1(a)=\bigl\{ k\in \{1,\dots,\mu_a\}| \frac{2\pi k\nu_a}{\mu_a}\text{ } 
mod
\text{ }2\pi\in (0,\pi]\bigr\}$, 
$N_2(a)=\bigl\{ k\in \{1,\dots,\mu_a\}| \frac{2\pi k\nu_a}{\mu_a}\text{ } mod
\text{ }2\pi\in (0,\pi)\bigr\}$. 

Define $R_a\in \Z[H_1(N)]$ by the
following formula:

\begin{multline}
R_a=
\sum_{k\in N_1(a)}\Bigl([K]f_a^{\mu_a-k}-f_a^k+
f_a^{k-\mu_a}-[K]f_a^{-k}\Bigr)\\
-\sum_{k\in N_2(a)}\Bigl(f_a^{\mu_a-k}-[K]f_a^k+
[K]f_a^{k-\mu_a}-f_a^{-k}\Bigr)
\end{multline}

Let $\tilde H$ be the quotient of
$\Z[H_1(N)]$ (viewed as a $\Z$-module) by the free
$\Z$-submodule generated by $\Bigl\{(e-f)([K]+e),\{R_a\}_{a\in A}\Bigr\}$.

One can prove that under the change of the orientation of $\pi^{-1}(a)$ 
(used to define $f_a$) $R_a$ goes to 
$-R_a$. Thus $\tilde H$ is well defined. 
To show this, one verifies that if $\mu_a$ is odd, then $N_1(a)=N_2(a)$. Under
this change each term from the first sum (used to define $R_a$) goes to minus
the corresponding term from the second sum and vice versa. (Note that
$f^2=e$.) If $\mu_a=2l$ is even, then $N_1(a)\setminus\{l\}=N_2(a)$. 
Under this change each term with $k\in N_1(a)\setminus \{l\}$ goes to minus the
corresponding term with $k\in N_2(a)$ and vice versa. 
The term in the first sum
that corresponds to $k=l$ goes to minus itself.

Finally define $\tilde S_K\in \tilde H$ as 
the sum over all double points $v$ of $\pi(K)$:
\begin{equation}\tilde S_K=
\sum_v\Bigl([\mu^+_1(v)]-[\mu^-_1(v)]+
[\mu^+_2(v)]-[\mu^-_2(v)]
\Bigr).
\end{equation}   
\end{defin}

\begin{thm}\label{correcttildeSSK}
$\tilde S_K$ is an isotopy invariant of $K$.
\end{thm}

The proof is a straightforward  generalization of the proof of
Theorem~\ref{correctSSK}.

%\begin{emf}\label{homversSSK}
%There is a version of $S_K$ with values in the free $\Z$-module
%generated by all free homotopy classes of oriented curves in $N$.
%To define it, one replaces all 
%homology classes in the definitions with the corresponding free homotopy
%classes. The summation should be made 
%over all double points of $\pi(K)$ such that
%neither of the knots obtained by the splitting is homotopic to a
%trivial loop and neither one of them is homotopic to a 
%positively oriented fiber of $N$.
%Each homology class of the form 
%$f_a^k$ in the definitions of $R^1_a$ and $R^2_a$ should be
%substituted with the free homotopy class of the curve that goes along the
%exceptional fiber $k$  times. Each class of the form $[K]f_a^k$ should be
%substituted with the free homotopy class of $K$ with a curve 
%going along the exceptional fiber $k$ times added to it. 
%(Note that $f_a$ is in the center
%of $\pi_1(N)$, so this class is well defined.)
%To prove that this is indeed an invariant of $K$, one easily modifies the
%proof of Theorem~\ref{correctSSK}.

%In a similar way one constructs
%the version of $\tilde S_K$ with values in the free $\Z$-module 
%generated by the set of all free homotopy classes of oriented 
%curves in $N$.
%$\end{emf} 

\begin{emf} 
One can easily verify that $S_K$ and $\tilde S_K$ 
satisfy relations similar to~\eqref{type2}. 
Hence both of them are Vassiliev invariants of order one
(see~\ref{vassiliev2}). 

%The corresponding versions of
%Theorem~\ref{realization2} also hold.
\end{emf}
 
\section{Wave fronts on surfaces}
\subsection{Definitions}
Let $F$ be a two-dimensional manifold. A {\em contact element\/}
at a point in $F$ is a one-dimensional vector subspace of the tangent plane. 
This subspace divides the tangent plane into two half-planes. A choice of
one of them is called a {\em coorientation\/} of a contact element. 
The space of all cooriented contact
elements of $F$ is a spherical cotangent bundle $ST^*F$. 
We will also denote it by $N$. It is an $S^1$-fibration over $F$. The
natural contact structure on $ST^*F$ is a
distribution of hyperplanes given by the condition that a velocity vector of
an incidence point of a contact element belongs to the element. 
A {\em Legendrian\/} curve $\lambda$ in $N$ is an immersion of $S^1$ into
$N$ such that for each $p\in S^1$ 
the velocity vector of $\lambda$ at $\lambda(p)$ 
lies in the contact plane. The naturally cooriented 
projection $L\subset F$ of a Legendrian curve $\lambda\subset N$ 
is called {\em the wave front\/} of
$\lambda$. 
A cooriented wave front may be uniquely lifted to a Legendrian curve
$\lambda \subset
N$ by taking a coorienting normal direction as a contact element at each
point of the front.
A wave front is said to be generic if it is an immersion
everywhere except a finite number of points, where it has cusp
singularities,
and all multiple points are double points with transversal self-intersection.
A cusp is the projection of a point where the corresponding Legendrian
curve is tangent to the fiber of the bundle.

\subsection{Shadows of wave fronts}

\begin{emf}\label{orientST*F}
For any surface $F$ the space $ST^*F$ is canonically oriented. The
orientation is constructed as follows. 
For a point $x\in F$ fix an
orientation of $T_xF$. It induces an orientation of the fiber over $x$. 
These two orientations determine an orientation of three dimensional planes
tangent to the points of the fiber over $x$. A straightforward 
verification shows 
that this orientation is independent on the orientation of $T_xF$ we choose.
Hence the orientation of $ST^*F$ is well defined.

Thus for oriented $F$ the shadow of a generic knot in $ST^*F$ is well defined 
(see~\ref{prelim} and~\ref{shadowgeneral}). Theorem~\ref{front-shadow}
describes the shadow of a Legendrian lifting of a generic cooriented wave front
$L\subset F$.  
\end{emf}

\begin{figure}[htb]
\begin{center}
 \epsfxsize\hsize\advance\epsfxsize -0.5cm
 \hepsffile{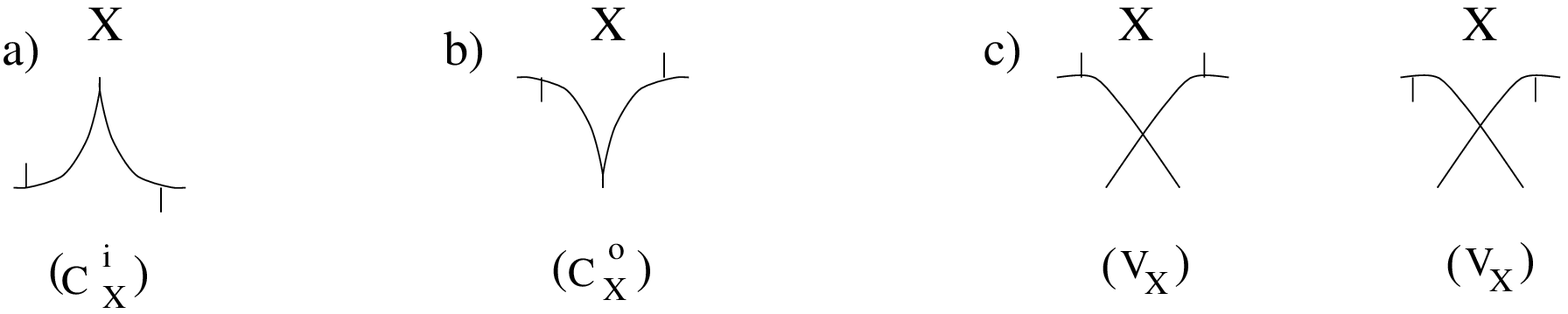}
 \end{center}
\caption{}\label{pic2.fig}
\end{figure}

\begin{defin} Let $X$ be a connected component of $F\setminus L$. We denote 
by $\chi \Int (X)$ the Euler characteristic of $\Int (X)$, 
by $C^i_X$ the number of cusps in the boundary of the region $X$ pointing 
inside
$X$ (as in Figure~\ref{pic2.fig}a), by $C^o_X$ the number of cusps 
in the boundary of $X$
pointing outside (as in Figure~\ref{pic2.fig}b),
and by $V_X$ the number of corners of $X$ where locally the picture
looks in one the two ways shown in Figure~\ref{pic2.fig}c. It
can happen that a cusp point is pointing both inside and outside of $X$. In
this case it contributes both in $C^i_X$ and in $C^o_X$. If the corner of the
type shown in Figure~\ref{pic2.fig}c enters twice in $\p X$, 
then it should be counted twice.  
\end{defin}

\begin{thm}\label{front-shadow}
Let $F$ be an oriented surface and
$L$ a generic cooriented wave front on $F$ corresponding to a 
Legendrian curve $\lambda$. 
There exists a small deformation of $\lambda$ in the class of all smooth 
(not 
only Legendrian) curves such that the resulting curve is generic with respect
to the projection, and the shadow of this curve can be constructed in the 
following way. We 
replace a small neighborhood of each cusp of $L$ with a smooth simple arc. 
The gleam of an arbitrary region $X$ that has a compact closure and does not
contain boundary components of $F$ is calculated by the following formula:
\begin{equation}
\gl_X=\chi \Int (X)+\frac {1}{2}(C^i_X-C^o_X-V_X).
\end{equation}
\end{thm}

For the proof of Theorem~\ref{front-shadow} see 
Subsection~\ref{pffront-shadow}.

\begin{rem}
The surface $F$ in the statement of Theorem~\ref{front-shadow} is not
assumed to be compact.

Note that as we mentioned in~\ref{shadowgeneral}, the gleam of a region $X$ 
that does not have compact closure or contains boundary components is
not well defined unless we fix a section over all ends of $X$ and components of
$\p F$ in $X$.

This theorem first appeared in~\cite{Tchernov}. A similar result was
independently obtained by Polyak~\cite{MPolyak}. 
\end{rem}

\begin{figure} [htbp]
\begin{center}
 \epsfxsize 12.5cm
 \hepsffile{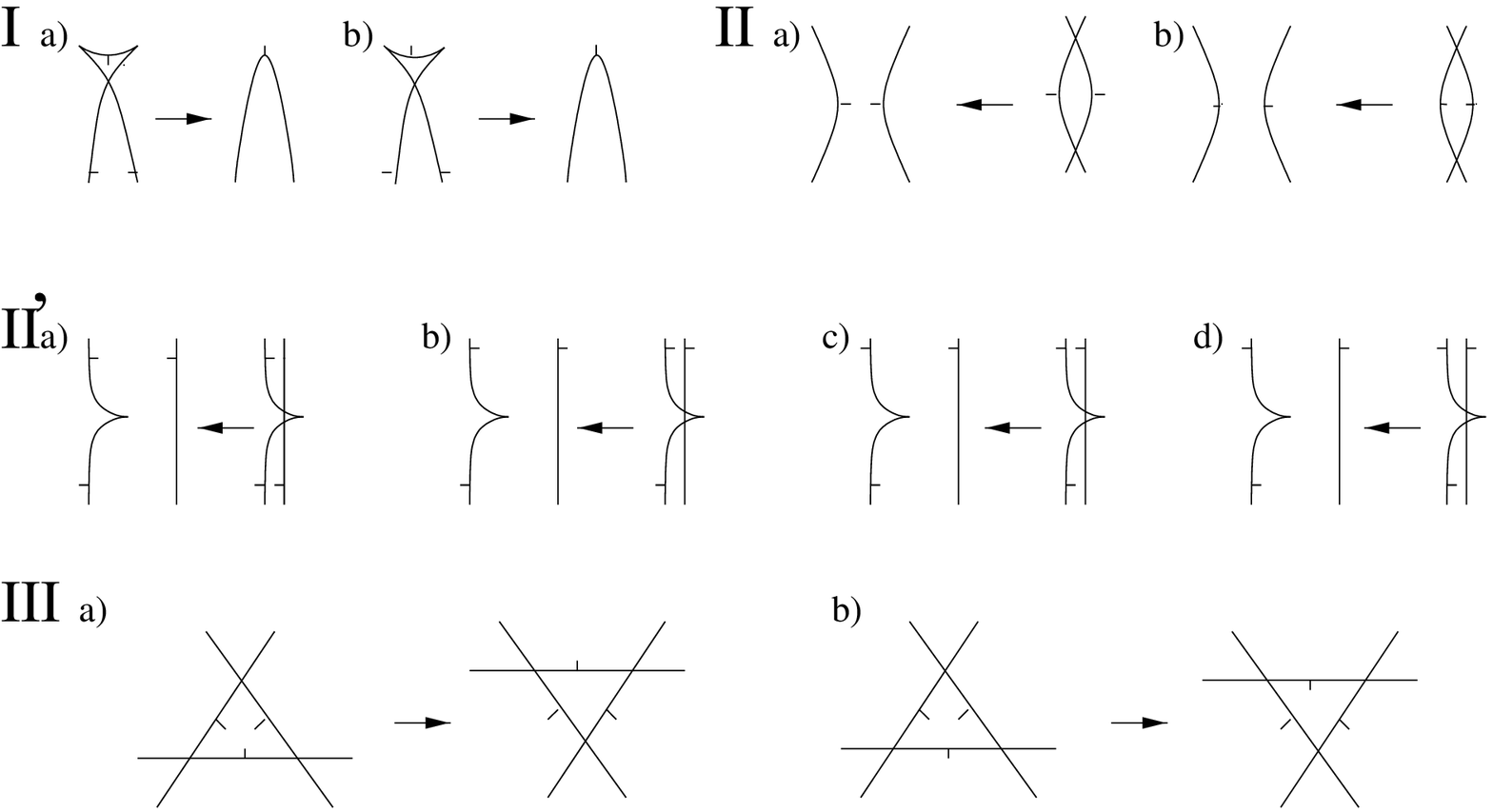}
 \end{center}
\caption{}\label{pic6.fig}
\end{figure}

\begin{emf} A self-tangency point $p$ of a wave front is said to be a
point of a {\em dangerous self-tangency\/} if the coorienting normals of the two 
branches coincide at $p$ (see Figure~\ref{pic8.fig}). Dangerous
self-tangency points correspond to self-intersection of the Legendrian curve. 
Hence a generic deformation of the front $L$ 
not involving {\em dangerous\/} self-tangencies
corresponds to an isotopy of the Legendrian knot $\lambda$.

\begin{figure}[htb]
\begin{center}
 \epsfysize 1.1cm
 \hepsffile{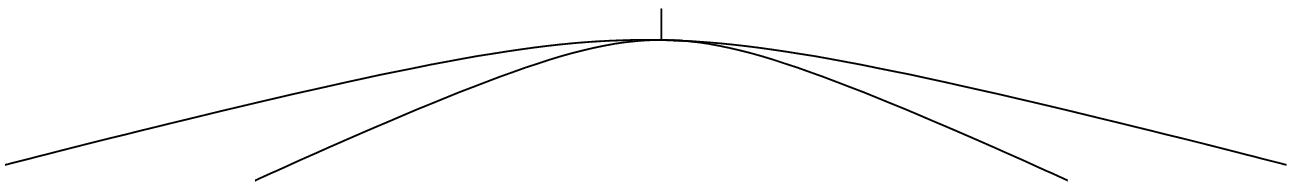}
 \end{center}
\caption{}\label{pic8.fig}
\end{figure}

Any generic deformation of a wave front $L$ corresponding to an isotopy
in the class of the Legendrian knots 
can be split into a sequence of modifications shown in 
Figure~\ref{pic6.fig}. The construction of Theorem~\ref{front-shadow}
transforms these generic modifications of wave fronts to shadow moves:
Ia and Ib in Figure~\ref{pic6.fig}  are transformed to the 
$\bar S_1$ move for shadow diagrams, IIa, IIb, II'a, II'b, II'c, and
II'd  are transformed to the $S_2$ move, 
finally IIIa and IIIb are transformed to $S_3$ and $\bar S_3$ respectively.
\end{emf}

\begin{emf} Thus for the Legendrian lifting of a wave front 
we are able to calculate all invariants that we can calculate for shadows. 
This includes the analogue of the linking number for the
fronts on $\R^2$ (see~\cite{Turaev}), the second order Vassiliev invariant
(see~\cite{Shumakovitch}), 
and quantum state sums (see~\cite{Turaev}).
\end{emf}

\subsection{Invariants of wave fronts on surfaces.}
In particular, the $S_K$ invariant gives rise to an invariant of a generic
wave front. This invariant appears to be related to the formula for the
Bennequin invariant of a wave front introduced by Polyak in~\cite{Polyak}.

Let us recall the corresponding results and definitions of~\cite{Polyak}.  

Let $L$ be a generic cooriented oriented wave front on an oriented surface
$F$.
A branch of a wave front is said to be positive (resp.
negative) if the frame of coorienting and  orienting vectors defines
positive (resp. negative) orientation of the surface $F$. 
Define the {\em sign\/} $\epsilon (v)$ of a double point $v$ of $L$ to be
$+1$ if the signs of both branches of the front intersecting at $v$
coincide and
$-1$ otherwise.  Similarly 
we assign a positive (resp. negative) sign to a cusp point
if the coorienting vector turns in a positive (resp. negative)
direction while traversing a small neighborhood of the cusp point along the
orientation. We denote half of the number of 
positive and negative cusp points by 
$C^+$ and $C^-$ respectively.

\begin{figure}[htb]
\begin{center}
 \epsfxsize 10cm 
 \hepsffile{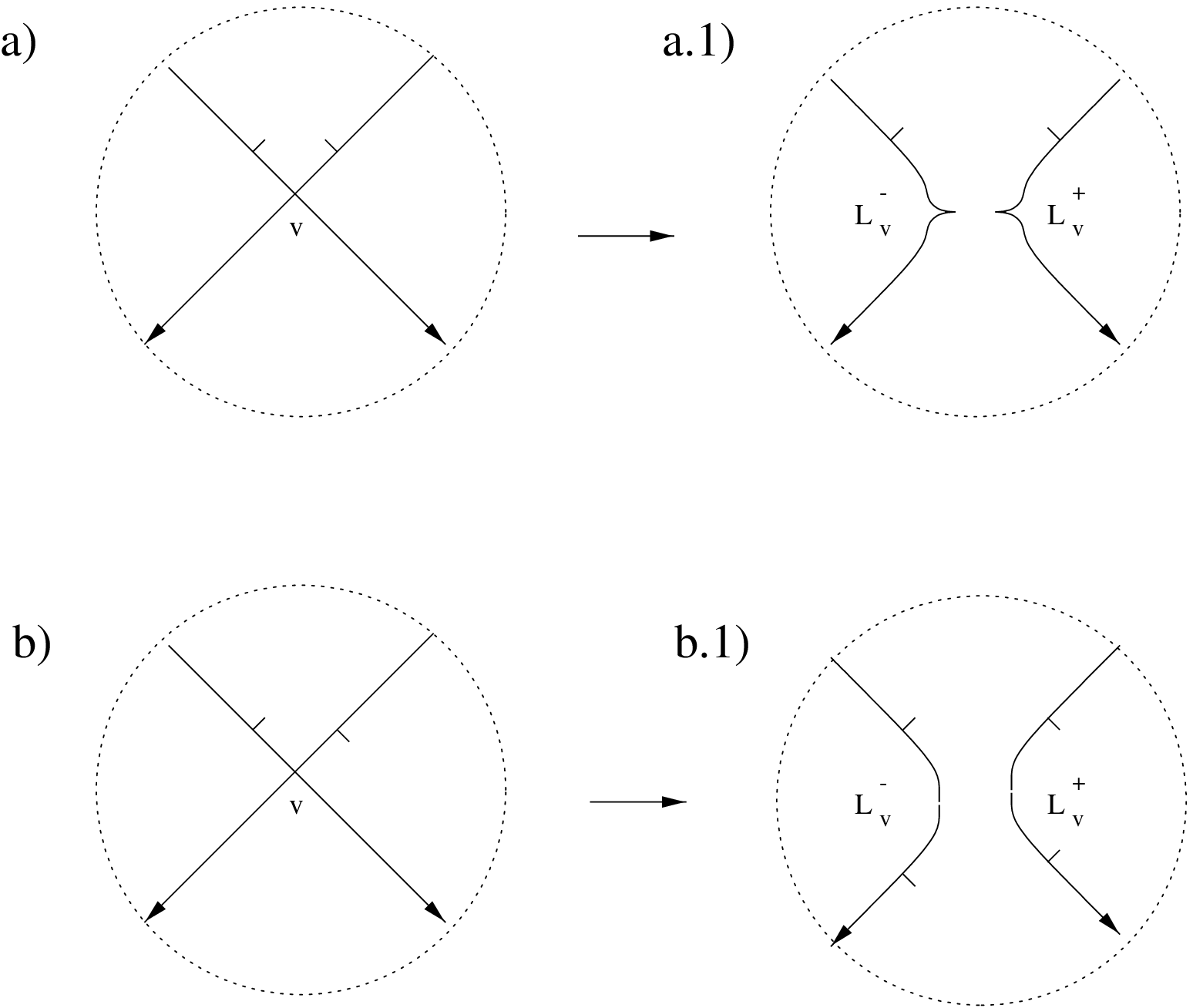}
 \end{center}
\caption{}\label{pic9.fig}
\end{figure}

Let $v$ be a double point of $L$. The orientations of $F$ and $L$ allow one
to distinguish the two wave fronts $L^+_v$ and $L^-_v$ obtained by splitting 
of $L$ in $v$ with respect to orientation and
coorientation (see Figures~\ref{pic9.fig}a.1 and~\ref{pic9.fig}b.1).
(Locally one of the two fronts lies to the left and another to the right of
$v$.) 

For a Legendrian curve $\lambda$ in $ST^*\R ^2$ denote by $l(\lambda)$ its
Bennequin invariant described in the works of 
Tabachnikov~\cite{Tabachnikov} and Arnold~\cite{Arnold} with the sign
convention of~\cite{Arnold} and~\cite{Polyak}. 

\begin{thm}[Polyak~\cite{Polyak}]\label{BenR2}
Let $L$ be a generic oriented cooriented wave front on $\R^2$
and $\lambda$ the corresponding Legendrian curve. Denote by $\ind(L)$ 
the degree of the mapping taking a point $p\in S^1$   
to the point of $S^1$ corresponding to the direction of 
the coorienting normal of $L$ at $L(p)$. Define $S$ as the 
following sum over all double points of $L:$
\begin{equation}
S=\sum_v (\ind(L^+_v)-\ind(L^-_v)-\epsilon (v)).
\end{equation}
Then 
\begin{equation}
l(\lambda )=S+(1-\ind(L))C^++(\ind(L)+1)C^- +\ind(L)^2.
\end{equation}
\end{thm}

In~\cite{Polyak} it is shown that the Bennequin invariant of a wave front on
the $\R^2$ plane admits quantization. Consider a formal quantum parameter
$q$. Recall that for any $n\in \Z$ the corresponding quantum number
$[n]_q\in \Z[q,q^{-1}]$ is a finite Laurent polynomial in $q$ defined by 
\begin{equation}
[n]_q=\frac{q^n-q^{-n}}{q-q^{-1}}.
\end{equation}
Substituting  quantum integers instead of integers in~\ref{BenR2}
we get the following theorem.

\begin{thm}[Polyak~\cite{Polyak}]\label{QBen}
Let $L$ be a generic cooriented oriented wave 
front on $\R^2$ and $\lambda$ the corresponding
Legendrian curve. Define $S_q$ by the following
formula, where the sum is taken over the set of all 
double points of $L:$
\begin{equation}\label{Sq}
S_q=\sum_v[ind(L^+_v)-ind(L^-_v)-\epsilon (v)]_q.
\end{equation}
Put
\begin{equation}
l_q(L)=S_q+[1-ind(L)]_qC^++[ind(L)+1]_qC^-+[ind(L)]_qind(L).
\end{equation} 
Then $l_q(\lambda)=l_q(L)\in\frac{1}{2}\Z[q,q^{-1}]$ is invariant 
under isotopy in the class of the Legendrian knots. 
\end{thm}

The $l_q(\lambda)$ invariant can be expressed~\cite{Aicpriv} 
in terms of the partial linking
polynomial of a generic cooriented oriented wave front
introduced by Aicardi~\cite{Aicardi}.

The reason why this invariant takes values in $\frac{1}{2}\Z[q,q^{-1}]$ and
not in $\Z[q,q^{-1}]$ is that the number of positive (or negative) cusps can
be odd. This makes $C^+$ ($C^-$) a half-integer.

Let $\lambda ^{\epsilon}_v$ with $\epsilon =\pm$ be the Legendrian lifting 
of the front $L^{\epsilon}_v$. Let $f\in 
H_1(ST^*F)$ be the homology class of a positively
oriented fiber.

\begin{thm}[Polyak~\cite{Polyak}]\label{Ben} Let $L$ be a generic oriented 
cooriented wave front on an oriented surface $F$. Let $\lambda$ be the corresponding 
Legendrian curve. Define $l_F(\lambda)\in H_1(ST^*F,\frac{1}{2}\Z)$ 
by the following formula: 
\begin{equation}
l_{F}(\lambda)=\Bigl(\prod_v[\lambda ^+_v][\lambda ^-_v]^{-1}
f^{-\epsilon (v)}\Bigr)(f[\lambda]^{-1})^{C^+}([\lambda]f)^{C^-}
\end{equation}
(We use the multiplicative notation for the group operation in
$H_1(ST^*F)$.)

Then $l_{F}(\lambda)$ is invariant under isotopy in the class of the
Legendrian knots.
\end{thm}

The proof is straightforward. One verifies that $l_{F}(\lambda)$ is 
invariant under all oriented versions of non-dangerous self-tangency, 
triple point, cusp crossing, and cusp birth moves of the wave front.

In~\cite{Polyak} this invariant is denoted by $I_\Sigma^+(\lambda)$ and, in
a sense, it
appears to be a natural generalization of Arnold's $J^+$
invariant~\cite{Arnold} 
to the case of an
oriented cooriented wave front on an oriented surface.

Note that in the situation of Theorem~\ref{BenR2} the indices of all
the fronts involved are the images of the homology classes of their
Legendrian liftings under the natural identification of $H_1(ST^*\R
^2)$ with $\Z$. If one replaces everywhere in~\ref{BenR2} indices with 
the corresponding homology classes and puts $f$ instead of $1$, then the
only difference between the two formulas is the term $\ind ^2(L)$. (One has
to remember that we use the multiplicative notation for the operation 
in $H_1(ST^*F)$.)
 
\begin{figure}[htb]
\begin{center}
 \epsfxsize 12 cm
 \hepsffile{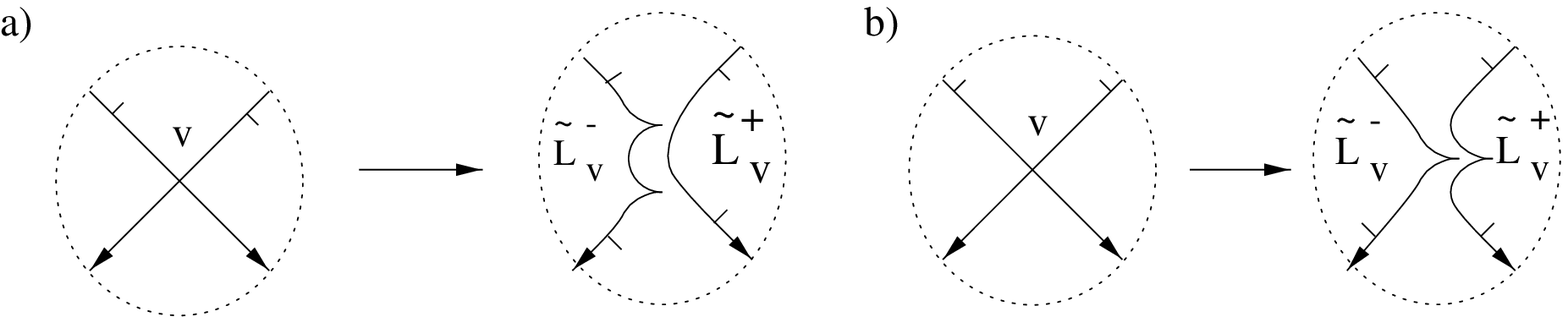}
 \end{center}
\caption{}\label{pic5.fig}
\end{figure}

\begin{emf} The splitting of a Legendrian knot $K$ into 
$\mu^+_v$  and $\mu^-_v$
(see~\ref{SK}) can be done up to an isotopy in the
class of the Legendrian knots. Although this can be done in many ways, there 
exists the simplest way. The projections $\tilde L^+_v$ and $\tilde L^-_v$
of the Legendrian curves created by the splitting are shown in 
Figure~\ref{pic5.fig}. (This fact follows from Theorem~\ref{front-shadow}.) 

Let $\tilde\lambda ^{\epsilon}_v$ with $\epsilon =\pm$ be the Legendrian lifting 
of the front $\tilde L^{\epsilon}_v$.
\end{emf}

\begin{thm}\label{splitwave} Let $L$ be a generic oriented cooriented wave
front on an oriented surface $F$. Let $\lambda$ be the corresponding
Legendrian curve. 
Define $S(\lambda)\in \frac{1}{2}\Z\bigl[ H_1(ST^*F)\bigr]$ by the following 
formula: 
\begin{equation}\label{SFL}
S(\lambda)=\sum_v\Bigl([\tilde \lambda ^+_v]-[\tilde \lambda
^-_v]\Bigr)
+(f-[\lambda])C^+ +([\lambda]f-e)C^-.
\end{equation} 
Then $S(\lambda)$ is invariant under isotopy in the class of the Legendrian
knots. 
\end{thm}

The proof is straightforward. One verifies that $S(\lambda)$ is indeed
invariant under all oriented versions of non-dangerous self-tangency, 
triple point, cusp crossing,  and cusp birth moves of the wave front.

\begin{emf}\label{homvers}
By taking the free homotopy classes of $\tilde \lambda^+_v$ and $\tilde
\lambda^-_v$ instead of the homology classes one obtains a different 
version of the $S(\lambda)$ invariant. It takes values in the group of
formal half-integer linear combinations of 
the free homotopy classes of oriented 
curves in $ST^*F$.
In this case the terms $[\lambda]$ and $f$ in~\eqref{SFL} should be
substituted with the free homotopy classes of $\lambda$ and of a 
positively oriented fiber respectively. The terms $[\lambda]f$ and $e$
in~\eqref{SFL} should be substituted with the free homotopy classes of 
$\lambda$ with a positive fiber added to it and the class of a 
contractible curve
respectively. Note that $f$ lies in the center of $\pi_1(ST^*F)$, so that the
class of $\lambda$ with a fiber added to it is well defined.

A straightforward verification shows that this version of $S(\lambda)$ is also
invariant under isotopy in the class of the Legendrian knots.
\end{emf}

\begin{thm}\label{splitting} Let $L$ be a generic oriented cooriented wave
front on an oriented surface $F$. Let $\lambda$ be the corresponding 
Legendrian curve. 
Let $S(\lambda)$ and $l_F(\lambda)$ be the invariants
introduced in~\ref{splitwave} and~\ref{Ben} respectively. Let
\begin{equation}
\pr:\frac{1}{2}\Z\bigl [H_1(ST^*F)\bigr]\rightarrow 
H_1(ST^* F,\frac{1}{2}\Z)
\end{equation}
be the  mapping defined as
follows: for any $n_i\in \frac{1}{2}\Z$ and $g_i\in H_1(ST^*F),$
\begin{equation}
\sum n_i g_i\mapsto \prod g_i^{n_i}.
\end{equation}
Then $\pr(S(\lambda))=l_F(\lambda)$.
\end{thm}
  
The proof is straightforward: one must verify that 
\begin{equation}\label{svyaz}
[\lambda ^+_v][\lambda ^-_v]^{-1}f^{-\epsilon (v)}=
[\tilde \lambda ^+_v][\tilde \lambda ^-_v]^{-1}
\text{ in } H_1(ST^*F). 
\end{equation}
(Recall that we use a multiplicative notation for the group operation
in $H_1(ST^*F)$.)

This means that $S_F(\lambda)$ is a refinement of Polyak's invariant 
$l_F(\lambda)$.

\begin{emf}\label{quant} 
One can verify that there is a unique linear combination 
$\sum_{m\in \Z}n_m[m]_q=l_q(\lambda)$ with $n_m$ being non-negative
half-integers such that $n_0=0$, and if $n_m>0$, then $n_{-m}=0$. 
To prove this one must verify that
$\{\frac{1}{2}[n]_q|0< n\}$ is a basis for the $\Z$-submodule of 
$\frac{1}{2}\Z[q,q^{-1}]$ 
generated by the quantum numbers and use the identity $n[m]_q=-n[-m]_q$. 
\end{emf}

The following theorem shows that if $L\subset \R^2$, then $S(\lambda)$ and
Polyak's quantization  $l_q(\lambda)$ (see~\ref{QBen}) of the Bennequin 
invariant can be explicitly expressed in terms of each other. 

\begin{thm}\label{equivalence}
Let $f\in H_1(ST^*\R^2)$ be the class of a positively oriented fiber. Let 
$L$ be a generic oriented cooriented wave front on $\R^2$, 
$\lambda$ the corresponding Legendrian
curve, and $f^h$ the homology class realized by it. Let
$l_q(\lambda)-[h]_qh=\sum_{m\in Z} n_m[m]_q$ be written in the form described 
in~\ref{quant} and $S(\lambda)=\sum_{l\in Z}k_l f^l$. 
Then
\begin{equation}\label{rel1}
l_q(\lambda)=[h]_qh+\sum_{k_l>0} k_l[2l-h-1]_q,
\end{equation}
and
\begin{equation}\label{rel2}
S(\lambda)=\sum_{n_m > 0}
n_m(f^{\frac{h+1+m}{2}}-f^{\frac{h+1-m}{2}}).
\end{equation}
\end{thm}
For the proof of Theorem~\ref{equivalence} see Subsection~\ref{pfequivalence}.

One can show that for $n_m>0$ both $\frac{h+1+m}{2}$ and $\frac{h+1-m}{2}$
are integers, so that the sum~\eqref{rel2} takes values 
in $\frac{1}{2}\Z[H_1(ST^*\R^2)]$.

Note that the $l_q(\lambda)$ invariant was defined only for fronts on the
plane $\R^2$. Thus $S(\lambda)$ is, in a sense, a generalization of
Polyak's
$l_q(\lambda)$ to the case of wave fronts on an arbitrary oriented surface $F$.

\begin{emf}
The splitting of the Legendrian knot $K$ into $\mu_1^+(v)$, $\mu_1^-(v)$, 
$\mu_2^+(v)$, and $\mu_2^-(v)$ (which was used to define $\tilde S(K)$, 
see~\ref{tildeSK}) can be done up to an isotopy in the class of the
Legendrian knots. Although this can be done in many ways, there is the
simplest one. The projections $\tilde L_1^+(v)$, $\tilde L_1^-(v)$, 
$\tilde L_2^+(v)$, and $\tilde L_2^-(v)$ are shown in
Figure~\ref{front1.fig}. (This fact follows from 
Theorem~\ref{front-shadow}.)

This allows us to introduce an invariant similar to $S(\lambda)$ 
for generic oriented cooriented 
wave fronts on a non-orientable surface $F$ in the
following way.

Let $L$ be a generic wave front on a
non-orientable surface $F$. Let $v$ be a double point of $L$. Fix some 
orientation of a small neighborhood of $v$ in $F$. The orientations of the
neighborhood and $L$ allow one to distinguish the wave fronts $L^+_1$,
$L^-_1$, $L^+_2$, and $L^-_2$ obtained by the two splittings of $L$ with
respect to the orientation and coorientation (see Figure~\ref{front1.fig}).
Locally the fronts with the upper indices plus and minus are located
respectively to the right and to the left of $v$.
To each double
point $v$ of $L$ we associate an element
$\Bigl([\tilde\lambda^+_1(v)]-[\tilde\lambda^-_1(v)]+[\tilde\lambda^+_2(v)]-
[\tilde\lambda^-_2(v)]\Bigr)\in \Z[H_1(ST^*(F)]$. Here we denote by lambdas
the Legendrian curves corresponding to the wave fronts appearing under
the splitting. Clearly 
this element does not depend on the orientation of the
neighborhood of $v$ we have chosen.

\begin{figure}[htb]
\begin{center}
 \epsfxsize 10cm
 \hepsffile{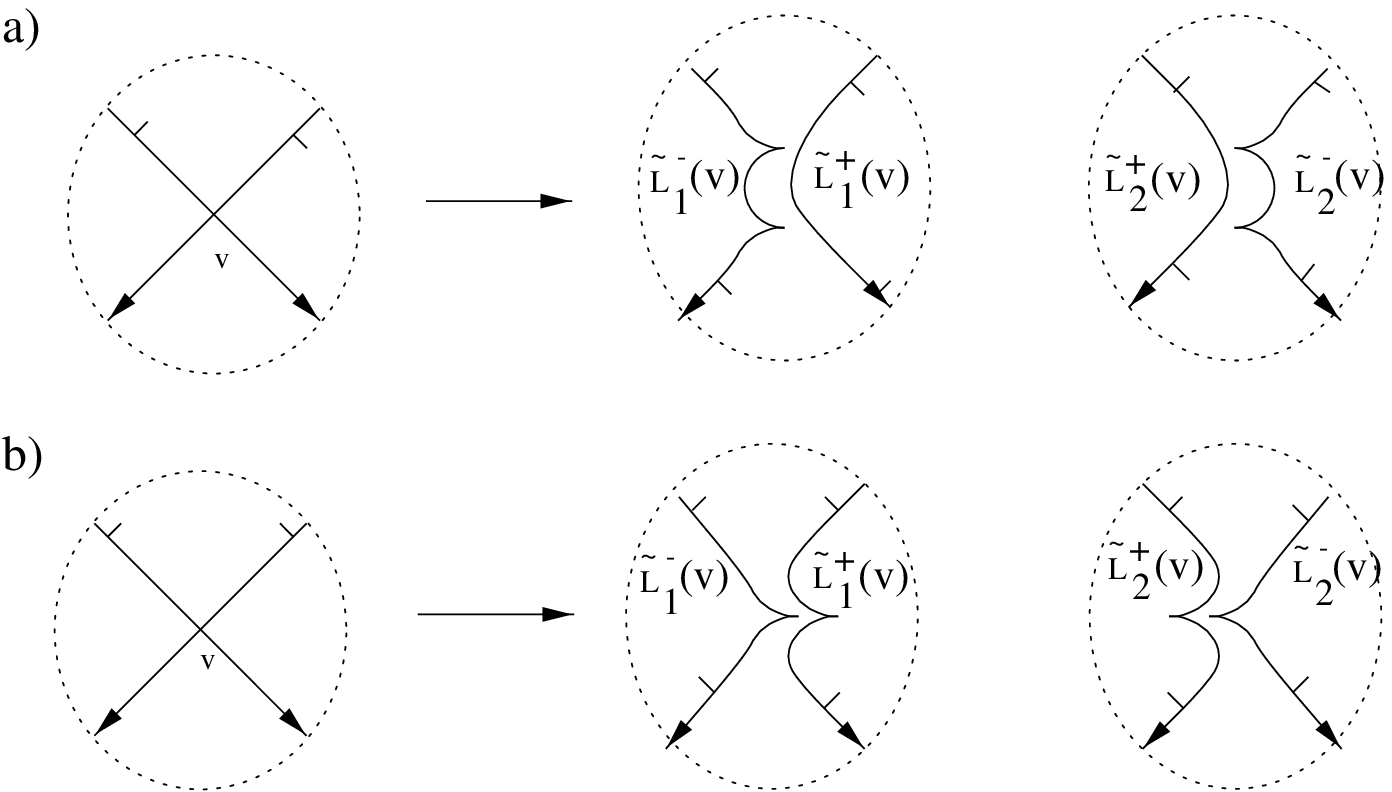}
 \end{center}
\caption{}\label{front1.fig}
\end{figure}

For a wave front $L$ let $C$ be half of the number of cusps of $L$. Denote
by $f$ the homology class of the fiber of $ST^*F$ oriented in some way. 
Note that $f^2=e$, so it does not matter which orientation of the fiber we
use to define $f$.
\end{emf}

\begin{thm}\label{unorientsplit}
Let $L$ be a generic cooriented oriented wave front on a non-orientable surface $F$
and $\lambda$ the corresponding Legendrian curve.
Define $\tilde S(\lambda)\in\frac{1}{2}\Z\bigl[ H_1(ST^*(F))\bigr]$ by the following
formula, where the summation is taken over the set of all double points of
$L:$
\begin{equation}\label{tildeSdef}
\tilde S(\lambda)=\sum_v\Bigl([\tilde \lambda ^+_1(v)]-
[\tilde \lambda ^-_1(v)]+[\tilde \lambda ^+_2(v)]-[\tilde \lambda
^-_2(v)]\Bigr)+C\bigl([\lambda]f-e+f-[\lambda]\bigr).
\end{equation}
Then $\tilde S(\lambda)$ is invariant under isotopy in the class of 
the Legendrian knots.
\end{thm}

The proof is straightforward. One verifies that $\tilde S(\lambda)$ is
indeed invariant under all oriented versions of non-dangerous
self-tangency, triple point passing, cusp crossing, and cusp birth moves of
the wave front.

The reason we have $\tilde S(\lambda)\in \frac{1}{2}\Z[H_1(ST^*F)]$ is that
if $L$ is an orientation reversing curve, then 
the number of cusps of $L$ is odd. In this case $C$ is a
half-integer.

%Similarly to~\ref{homvers}, one can introduce the version of the
%$\tilde S(\lambda)$ invariant with values in the group of formal
%half-integer linear combinations of free homotopy classes of
%oriented curves in $ST^*F$.

\section{Wave fronts on orbifolds}
\subsection{Definitions}
\begin{defin}\label{orbifold}

An {\em orbifold\/} is a surface $F$ with the additional structure
consisting of:  

1) set $A\subset F$;

2) smooth structure on $F\setminus A$;

3) set of homeomorphisms  $\phi_a$ of neighborhoods $U_a$ of $a$ in $F$ 
onto $\R^2/G_a$ such that $\phi_a(a)=0$ and
$\phi_a\Big|_{U_a\setminus a}$ is a diffeomorphism.
Here $G_a=\bigl\{e^{\frac{2\pi
k}{\mu_a}}\big| k\in\{1,\dots,\mu_a\}\bigr\}$ is a group
acting on $\R^2=\C$ by multiplication. ($\mu _a\in \Z$ is assumed to be
positive.) 
\end{defin}

The points $a\in A$ are called {\em cone points\/}.

The action of $G$ on $\R^2$ induces the action of $G$ on
$ST^*\R^2$. This makes $ST^*\R^2/G$ a Seifert fibration over
$\R^2/G$. Gluing together the pieces over neighborhoods of $F$ we obtain a
Seifert fibration $\pi:N\rightarrow F$. The fiber over a cone point $a$ is an
exceptional fiber of type $(\mu_a,-1)$ (see~\ref{numbers}).  
   
The natural contact structure on $ST^*\R^2$ is invariant under 
the induced action of $G$. Since $G$ acts freely on $ST^*\R^2$, this implies
that $N$ has an induced contact structure. As before, the naturally
cooriented projection
$L\subset F$ of a generic Legendrian curve $\lambda$ is called 
{\em the front of\/} $\lambda$.
  
\subsection{Invariants for fronts on orbifolds}
For oriented $F$ we construct an
invariant similar to $S(\lambda)$. It corresponds to the $S_K$
invariant of a knot in a Seifert fibered space. 
For non-orientable surface $F$
we construct an analogue of $\tilde S(\lambda)$. It corresponds to
the $\tilde S_K$ invariant of a knot in a Seifert fibered space.

Note that any surface $F$ can be viewed as an orbifold without cone
points. This justifies the notation below.

Let $F$ be an oriented surface. The orientation of $F$ induces an 
orientation of all fibers. Denote by $f$ 
the homology class of a positively oriented fiber. For a cone point $a$
denote by $f_a$ the homology class of a positively oriented fiber
$\pi^{-1}(a)$. For a generic oriented cooriented wave front $L\subset F$
denote by $C^+$ (resp. $C^-$) half of the number 
of positive (resp. negative) cusps of $L$.
Note that for a double point $v$ of a generic front $L$
the splitting into $\tilde L^+_v$ and $\tilde L^-_v$ is well defined. 
The corresponding Legendrian
curves $\tilde \lambda^+_v$ and $\tilde \lambda^-_v$ in 
$N$ are also well defined.

For $a\in A$ of type $(\mu_a,-1)$ 
put $N_1(a)=\bigl\{k\in\{1,\dots,\mu_a\}\big|\frac{-2k\pi}{\mu_a}
\mod 2\pi
\in(0,\pi]\bigr\}$, 
$N_2(a)=\bigl\{k\in\{1,\dots,\mu_a\}\big|\frac{2k\pi}{\mu_a}
\mod 2\pi
\in(0,\pi)\bigr\}$. 
Define $R^1_a, R^2_a\in
\Z[H_1(N)]$ by the following formulas: 

\begin{equation}
R^1_a=\sum_{k\in N_1(a)}\bigl([\lambda]f_a^{\mu_a-k}-f_a^k\bigr)-
\sum_{k\in N_2(a)}\bigl(f_a^{\mu_a-k}-[\lambda]f_a^k\bigr),
\end{equation}
\begin{equation}
R^2_a=\sum_{k\in N_1(a)}\bigl(f_a^{\mu_a-k}-[\lambda]f_a^k\bigr)-
\sum_{k\in N_2(a)}\bigl([\lambda]f_a^{\mu_a-k}-f_a^k\bigr).
\end{equation}

Set $J$ to be the quotient of
$\frac{1}{2}\Z[H_1(N)]$ by the free Abelian 
subgroup generated by 
$\Bigl\{\{\frac{1}{2}R_1(a),\frac{1}{2}R_2(a)\}_{a\in A}\Bigr\}$.

\begin{thm}\label{orbifold1}
Let $L$ be a generic cooriented oriented wave front on $F$ and
$\lambda$ the corresponding Legendrian curve.

Then $S(\lambda)\in J$ defined by the sum over 
all double points of $L$,
\begin{equation}
S(\lambda)=
\sum\Bigl([\tilde \lambda^+(v)]-[\tilde \lambda^-(v)]\Bigr)+(f-[\lambda])C^++
([\lambda]f-e)C^-,
\end{equation}
is invariant under isotopy in the class of the Legendrian knots.
\end{thm}  

For the proof of Theorem~\ref{orbifold1} see Subsection~\ref{pforbifold1}.

Let $F$ be a non-orientable surface. 
Denote by $f$ the homology class of a regular fiber oriented in some way. 
Note that $f^2=e$, so the orientation we use to define $f$ does not matter. 
For a cone point $a$ denote by $f_a$ 
the homology class of the fiber $\pi^{-1}(a)$ oriented in some
way. For a generic oriented cooriented wave front $L\subset F$
denote by $C$ half of the number of cusps of $L$.
Note that for a double point $v$ of a generic front $L$ 
the element 
$\bigl([\tilde\lambda^+_1(v)]-[\tilde\lambda^-_1(v)]+
[\tilde\lambda^+_2(v)]-[\tilde\lambda^-_2(v)]\bigr)\in \Z[H_1(N)]$ 
used to introduce
$\tilde S(\lambda)$ is well defined.
   
For $a\in A$ of type $(\mu_a,-1)$  
put $N_1(a)=\bigl\{k\in\{1,\dots,\mu_a\}\big |\frac{-2k\pi}{\mu_a}
\mod 2\pi
\in(0,\pi]\bigr\}$, and 
$N_2(a)=\bigl\{k\in\{1,\dots,\mu_a\}\big |\frac{-2k\pi}{\mu_a}
\mod 2\pi
\in(0,\pi)\bigr\}$. 
Define $R_a\in
\Z[H_1(N)]$ by the following formula:
\begin{multline}
R_a=
\sum_{k\in N_1(a)}\Bigl([\lambda]f_a^{\mu_a-k}-f_a^k+
f_a^{k-\mu_a}-[\lambda]f_a^{-k}\Bigr)\\
-\sum_{k\in N_2(a)}\Bigl(f_a^{\mu_a-k}-[\lambda]f_a^k+
[\lambda]f_a^{k-\mu_a}-f_a^{-k}\Bigr).
\end{multline}

Put $\tilde J$ to be the quotient of
$\frac{1}{2}\Z[H_1(N)]$ by a free Abelian subgroup generated by 
$\Bigl\{\{\frac{1}{2}R_a\}_{a\in A}\Bigr\}$.

Similarly to~\ref{tilde1SK}, one 
can prove that under the change of the orientation of $\pi^{-1}(a)$ 
(used to define $f_a$) $R_a$ goes to 
$-R_a$. Thus $\tilde J$ is well defined.

\begin{thm}\label{orbifold2}
Let $L$ be a generic cooriented oriented wave front on $F$ and
$\lambda$ the corresponding Legendrian curve.

Then $\tilde S(\lambda)\in \tilde J$ defined by 
the summation  over all double points of $L$,
\begin{equation}
\tilde S(\lambda)=
\sum\Bigl([\tilde \lambda^+_1(v)]-[\tilde \lambda^-_1(v)]+
[\tilde \lambda^+_2(v)]-[\lambda^-_2(v)]
\Bigr)+\bigl(([\tilde \lambda]f-e+f-[\lambda]\bigr)C,
\end{equation}
is invariant under isotopy in the class of the Legendrian knots.
\end{thm}  

The proof is a straightforward  generalization of the proof of
Theorem~\ref{orbifold1}.

%As before (see~\ref{homvers} and~\ref{homversSSK}) one can 
%introduce versions of $S(\lambda)$ and $\tilde S(\lambda)$
%with values in the 
%quotients of the group of all formal half-integer linear combinations of 
%free homotopy classes of oriented curves in $N$.      

\section{Proofs}

\subsection{Proof of Theorem~\ref{correct2}.}\label{pfcorrect2}
To prove the theorem it suffices to show that $S_K$ is invariant
under the elementary isotopies. They project to: 
a death of a double point, cancellation of two double points, and passing
through a triple point. 

To prove the invariance, we fix a homeomorphic to a closed disk 
part $P$ of $F$ containing the projection of one of the elementary
isotopies.
Fix a section over the boundary of $P$ such that the 
points of $K\cap\pi^{-1}(\p P)$ belong to the section. Inside $P$ we can
construct the Turaev shadow (see \ref{shadowgeneral}). 
The action of $H_1(\Int P)=e$  
on the set of isotopy types of links is trivial (see~\ref{action}). 
Thus the part of $K$ can be reconstructed in the unique way 
from the shadow over $P$. 
In particular, one can compare the homology classes
of the curves created by splitting at a double point inside $P$. 
Hence to prove the theorem, it suffices to verify the invariance under the
oriented versions of the moves $S_1,S_2$, and $S_3$.

\begin{figure}[htbp]
 \begin{center}
  \epsfxsize 6cm
  \hepsffile{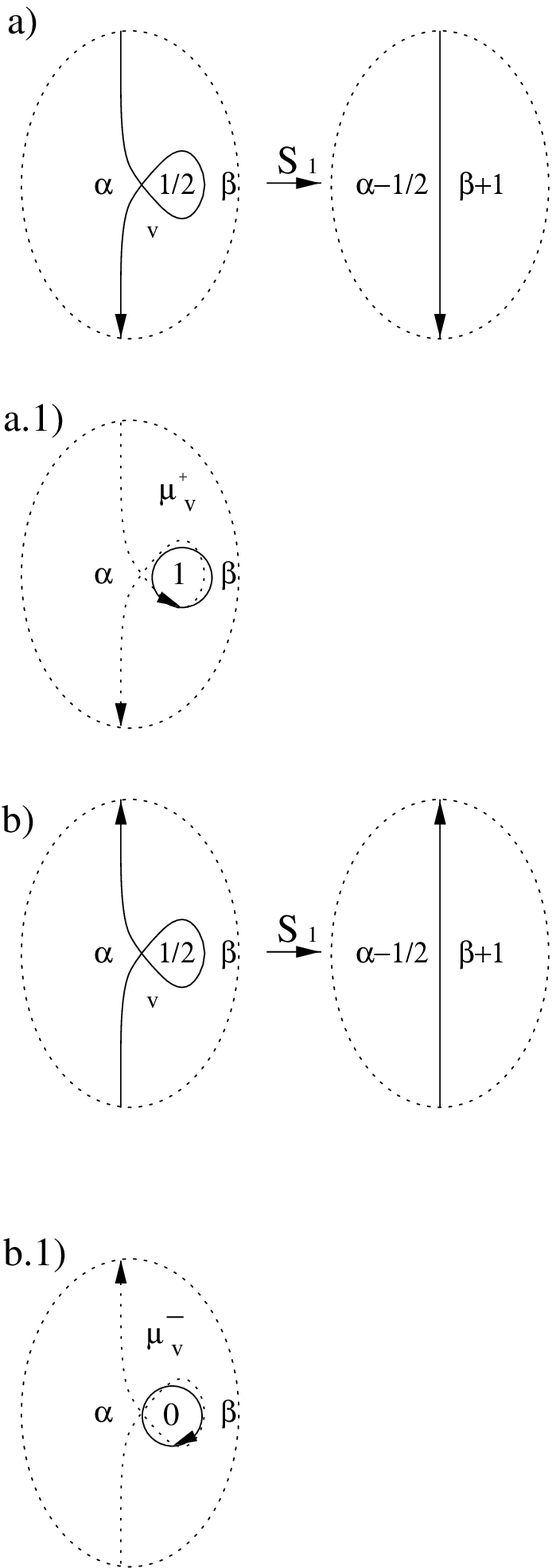}
 \end{center}
\caption{}\label{shad11.fig}
\end{figure}

There are two versions of the oriented move $S_1$ shown in
Figures~\ref{shad11.fig}a and~\ref{shad11.fig}b. 

For Figure~\ref{shad11.fig}a the term $[\mu^+_v]$ appears to be equal to
$f$. From~\ref{sum-property} we know that  
$[\mu^+_v][\mu^-_v]=[K]f$, so that $[\mu^-_v]=[K]$. Hence
$[\mu^+_v]-[\mu^-_v]=f-[K]$ and is equal to zero in $H$. In the same way we
verify that $[\mu^+_v]-[\mu^-_v]$ (for $v$ shown in
Figure~\ref{shad11.fig}b) is equal to $[K]f-e$. It is also zero in $H$. 
The summands corresponding to other double points do not change under this move, 
since it does
not change the homology classes of the knots created by the splittings.

\begin{figure}[htbp]
 \begin{center}
  \epsfxsize 8cm
  \hepsffile{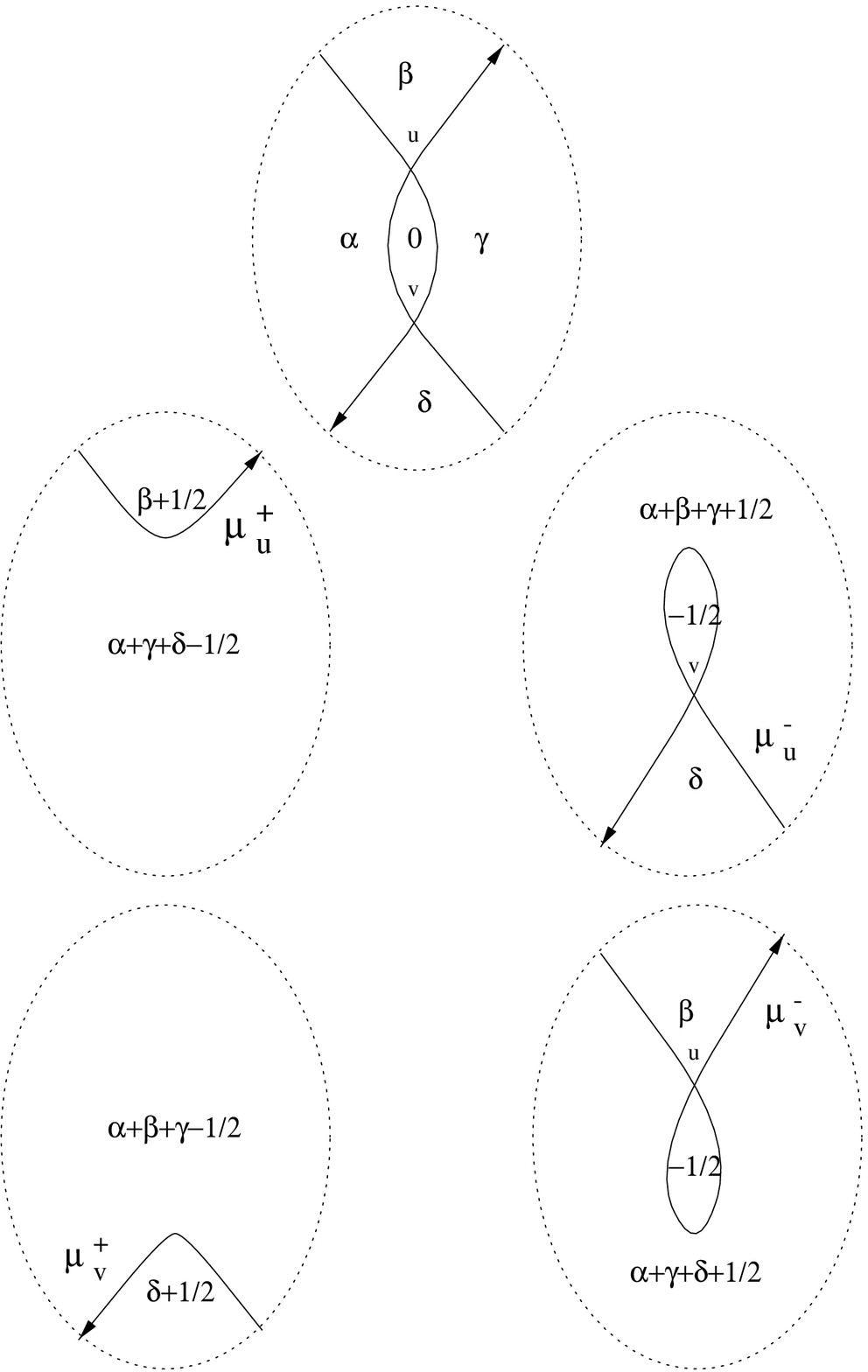}
 \end{center}
\caption{}\label{shad2.fig}
\end{figure}

   There are three oriented versions of the $S_2$ move. We show
that $S_K$ does not change under one of them. The proof for
the other two is the same or easier. We choose the version 
corresponding to the upper part of Figure~\ref{shad2.fig}. 
The summands corresponding to the double points not in this figure are 
preserved under the move, since it does not change the homology classes 
of the corresponding knots. So it suffices to show that the terms 
produced by 
the double points $u$ and $v$ in this figure cancel out. 
Note that the shadow $\mu ^-_v$ is transformed to $\mu^+_u$
by the $\bar S_1$ move. It is known that 
$\bar S_1$ can be expressed in terms of $S_1,S_2$, and 
$S_3$, thus it also does not change the homology classes of the 
knots created by the splittings. Hence $[\mu^+_u]$ and $[\mu^-_v]$ cancel out. 
In the same way one 
proves that $[\mu ^-_u]$ and $[\mu ^+_v]$ also cancel out. 

\begin{figure}[htb]
 \begin{center}
  \epsfxsize 7cm
  \hepsffile{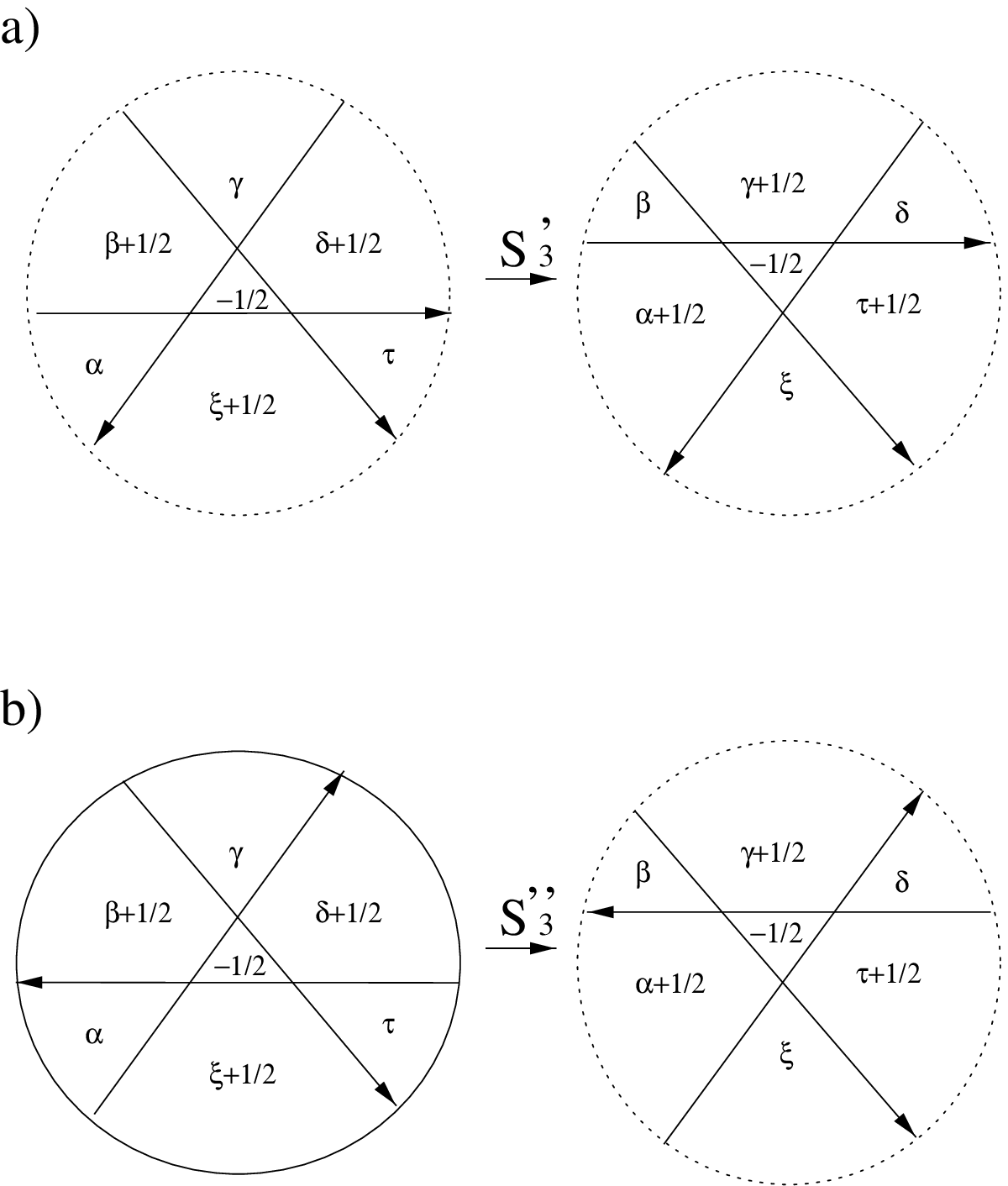}
 \end{center}
\caption{}\label{shad7.fig}
\end{figure}

\begin{figure}[htb]
 \begin{center}
  \epsfxsize 9cm
  \hepsffile{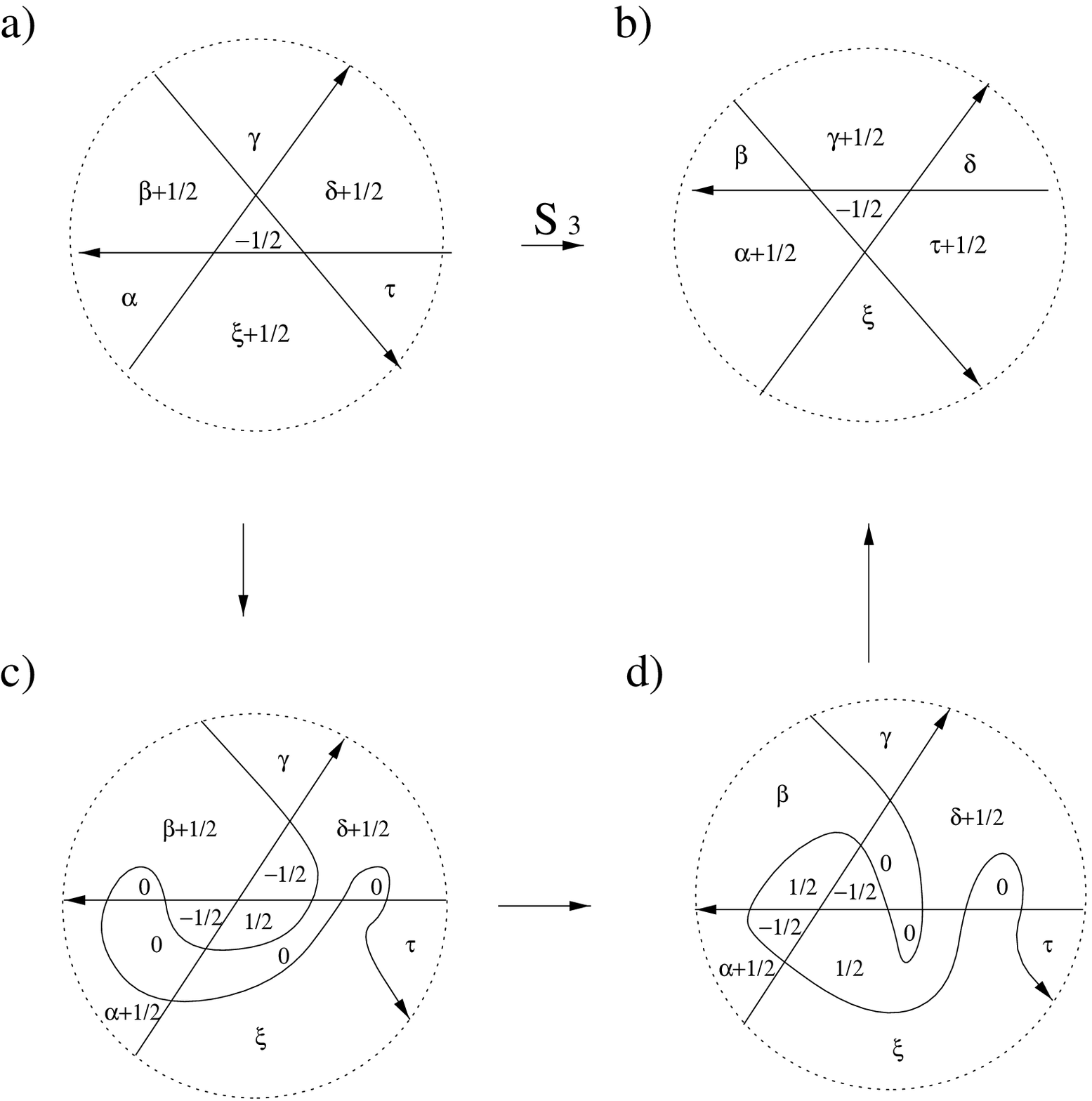}
 \end{center}
\caption{}\label{shad6.fig}
\end{figure}

There are two oriented versions of the $S_3$ move: $S'_3$ and 
$S''_3$, shown in Figures~\ref{shad7.fig}a and~\ref{shad7.fig}b respectively.
The $S''_3$ move can be expressed in terms of $S'_3$ and 
oriented versions of $S_2$ and $S_2^{-1}$.  To prove this we use 
Figure~\ref{shad6.fig}. There are two ways to get from 
Figure~\ref{shad6.fig}a to Figure~\ref{shad6.fig}b. One is to apply 
$S''_3$. Another way is to apply three times the oriented version of 
$S_2$ to obtain Figure~\ref{shad6.fig}c, 
then apply $S'_3$ to get Figure~\ref{shad6.fig}d, 
and finally use three times the oriented 
version of $S_2^{-1}$ to end up at Figure~\ref{shad6.fig}b.

\begin{figure}[htbp]
 \begin{center}
  \epsfxsize 10cm
  \hepsffile{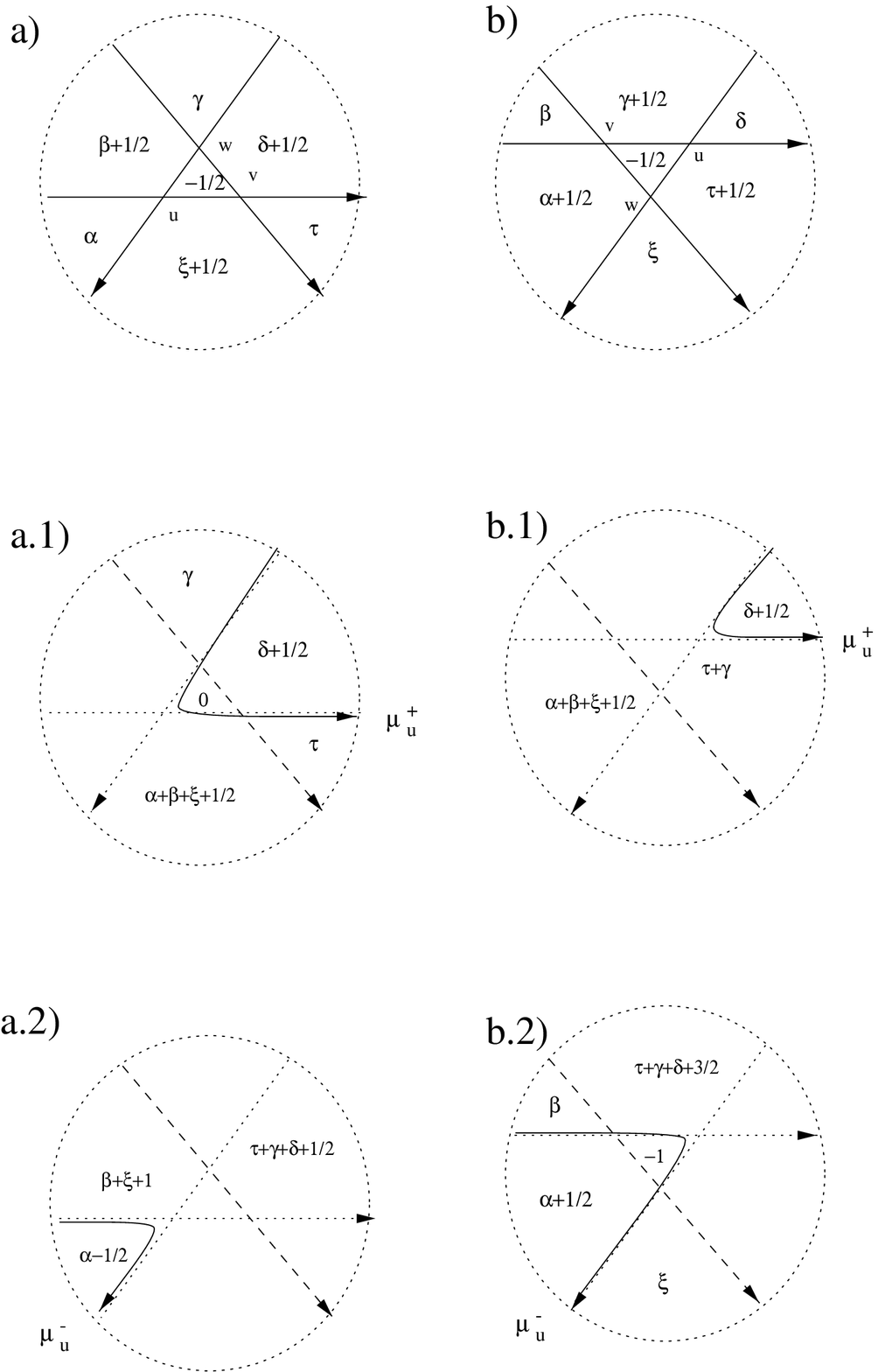}
 \end{center}
\caption{}\label{shad8.fig}
\end{figure}

Thus it suffices to verify the invariance under $S'_3$. The terms
corresponding to the double points not in Figures~\ref{shad8.fig}a 
and~\ref{shad8.fig}b are preserved for the same reasons as above. 
The terms coming from double points $u$ in 
Figure~\ref{shad8.fig}a and $u$ in Figure~\ref{shad8.fig}.b are the 
same. 
This holds also for the $v$- and $w$-pairs of double points in these two
figures. We prove this statement only for the $u$-pair of double points. For
$v$- and $w$-pairs the proof is the same or simpler. There is only
one possibility: either the dashed line belongs to both $\pi(\mu^+_u)$
in Figures~\ref{shad8.fig}a.1 and~\ref{shad8.fig}b.1 respectively or to
both $\pi(\mu ^-_u)$ in Figures~\ref{shad8.fig}a.2 
and~\ref{shad8.fig}b.2 respectively. 
We choose the one to which it does not belong. 
Summing up gleams on each of the two sides of it we immediately see
that the corresponding shadows are the same on both pictures. 
Thus the homology classes of the corresponding knots are equal. 
But $[\mu ^+_u][\mu ^-_u]=[K]f$ (see~\ref{sum-property}),
thus the homology classes of the knots represented by the 
other shadows are also equal. 
This completes the proof of Theorem~\ref{correct2}.\qed

\subsection{Proof of Theorem~\ref{realization2}.}\label{pfrealization2}
I: $K'$ can be obtained from $K$ by a sequence of isotopies and
modifications along fibers. Isotopies do not change $S$. The modifications 
change $S$ by
elements of type~\eqref{type2}. To complete the proof we use the
identity $\xi_1 \xi_2=[K]$.

II: We prove that for any $i\in H_1(N)$ there exist
two knots $K_1$ and $K_2$ such that they represent the same free homotopy
class as $K$, 
\begin{equation}
S_{K_1}=S_K+(f-e)([K]i^{-1}+i),\text{ and}
\end{equation}
\begin{equation} 
S_{K_2}=S_K-(f-e)([K]i^{-1}+i).
\end{equation}
Clearly this would imply the second statement of the theorem.

\begin{figure}[htbp]
 \begin{center}
  \epsfxsize 7cm
  \hepsffile{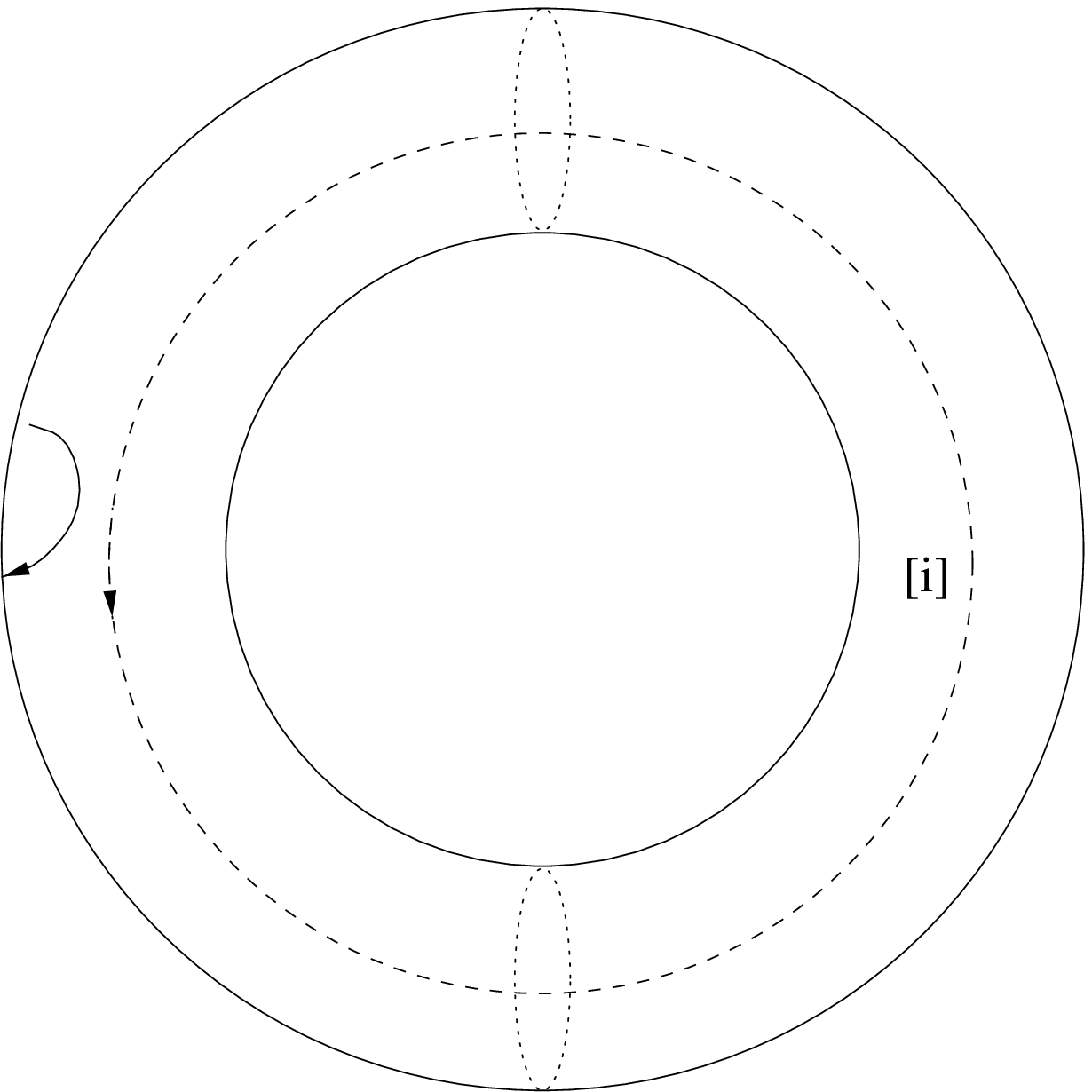}
 \end{center}
\caption{}\label{aicardi9.fig}
\end{figure}

\begin{figure}[htbp]
 \begin{center}
  \epsfxsize 7cm
  \hepsffile{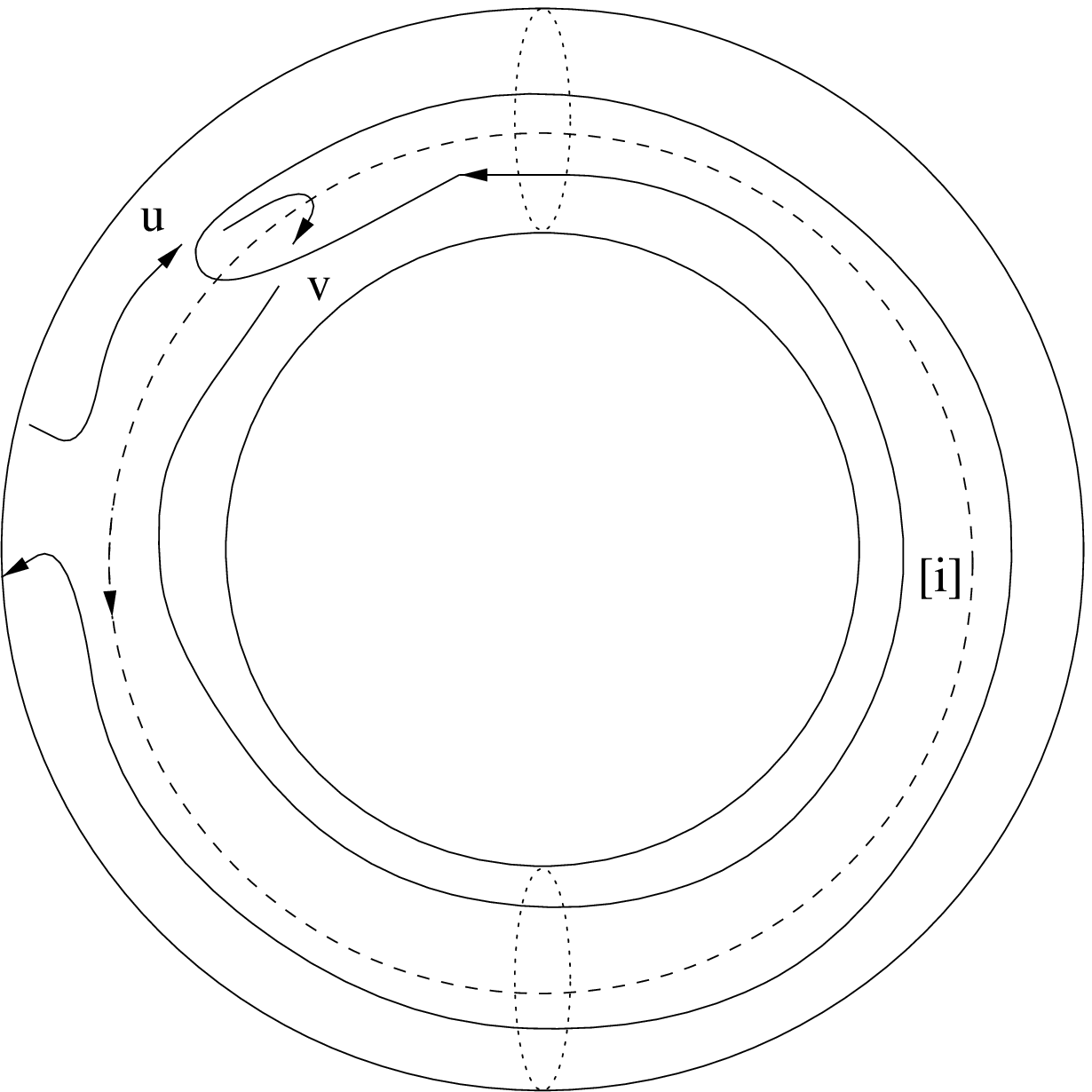}
 \end{center}
\caption{}\label{aicardi1.fig}
\end{figure}

Take $i\in H_1(N)$. Let $K_i$ be an oriented knot in $N$ such that $[K_i]=i$. 
The space $N$ is oriented, hence the tubular neighborhood $T_{K_i}$ of $K_i$ 
is homeomorphic to an oriented solid torus $T$. Deform $K_i$, so that
$K_i\cap T_{K_i}$ is a small arc (see Figure~\ref{aicardi9.fig}). 
Pull one part of the arc along $K_i$ in $T_{K_i}$ under the other part
of the arc (see Figure~\ref{aicardi1.fig}). This isotopy creates 
two new double points 
$u$ and $v$ of $\pi(K)$. (Since $T_{K_i}$ may be knotted, it might happen that 
there are other new double points, but we do not need them for
our construction.) Making a fiber modification along the part of 
$\pi^{-1}(u)$ that lies in $T$ one obtains $K_2$. Making a fiber
modification along the part of $\pi^{-1}(v)$ that lies in $T$ one obtains
$K_1$. 
This completes the proof of Theorem~\ref{realization2}.
\qed

\subsection{Proof of Theorem~\ref{homtor}.}\label{pfhomtor}
\begin{figure}[htbp]
 \begin{center}
  \epsfxsize 9cm
  \hepsffile{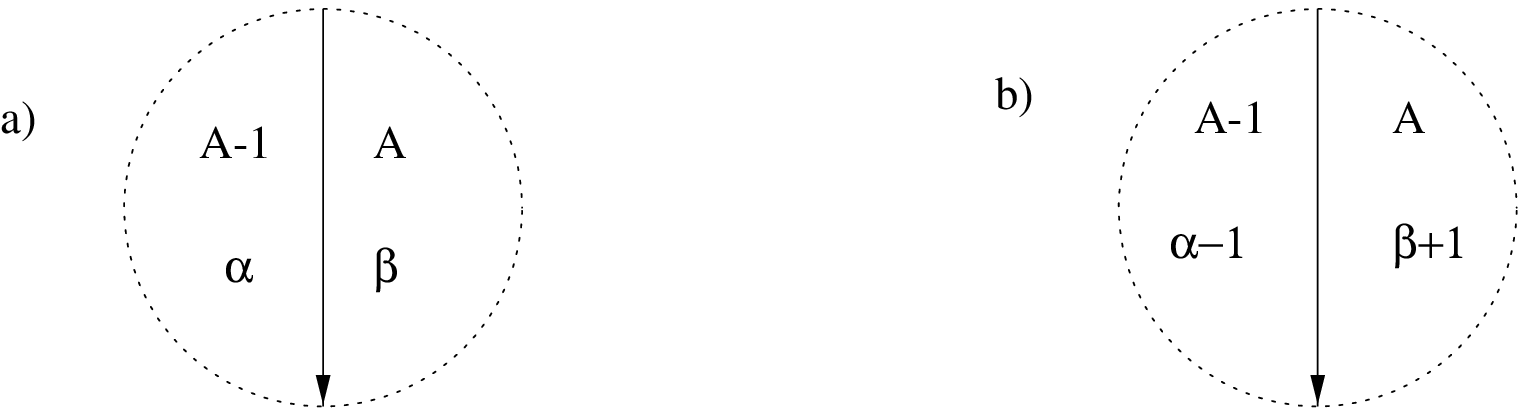}
 \end{center}
\caption{}\label{shad16.fig}
\end{figure}

It is easy to verify that any two shadows with the same projection can be 
transformed to each other by a sequence of fiber fusions. 
One can easily create a trivial knot with an ascending diagram 
such that its projection is any desired curve. 
This implies that any two shadows 
on $\R^2$ can be transformed to each other by a sequence of fiber fusions, 
movements $S_1,S_2,S_3$, and their inverses. 
A straightforward verification shows that $\sigma(s(K))$ does 
not change under the moves $S_1,S_2,S_3$, and their inverses. 
Under fiber fusions the homology class of the knot and the element 
$\sigma$ change in the same way. 
To prove this, we use Figure~\ref{shad16.fig}, where Figure~\ref{shad16.fig}a shows the shadow before the application 
of the fiber fusion (that adds $1$ 
to the homology class of the knot) and 
Figure~\ref{shad16.fig}b after. In this diagrams 
the indices of the regions are denoted by Latin letters. Now one easily 
verifies that $\sigma$ also increases by one. Finally, for 
the trivial knot with a trivial shadow diagram its homology class  
and $\sigma (s(K))$ are both equal to $0$.
This completes the proof of Theorem~\ref{homtor}.
\qed

\subsection{Proof of Theorem~\ref{correctSSK}.}\label{pfcorrectSSK}
It suffices to show that $S_K$ does
not change under the elementary isotopies of the knot. 
Three of them correspond in the projection to: a 
birth of a small loop, passing through a point of self-tangency, and passing
through a triple point.  The fourth one is 
passing through an exceptional fiber.

From the proof of Theorem~\ref{correct2} one gets that $S_K$
is invariant under the first three of the elementary isotopies described above. 
Thus it suffices to prove invariance under passing through an exceptional 
fiber $a$.

Let $a$ be a singular fiber of type $(\mu_a, \nu_a)$ (see~\ref{numbers}). 
Let $T_a$ be a neighborhood of $a$
which is fiber-wise isomorphic to the standardly fibered solid torus with
an exceptional fiber of type $(\mu_a, \nu_a)$. 

\begin{figure}[htbp]
 \begin{center}
  \epsfxsize 11cm
  \hepsffile{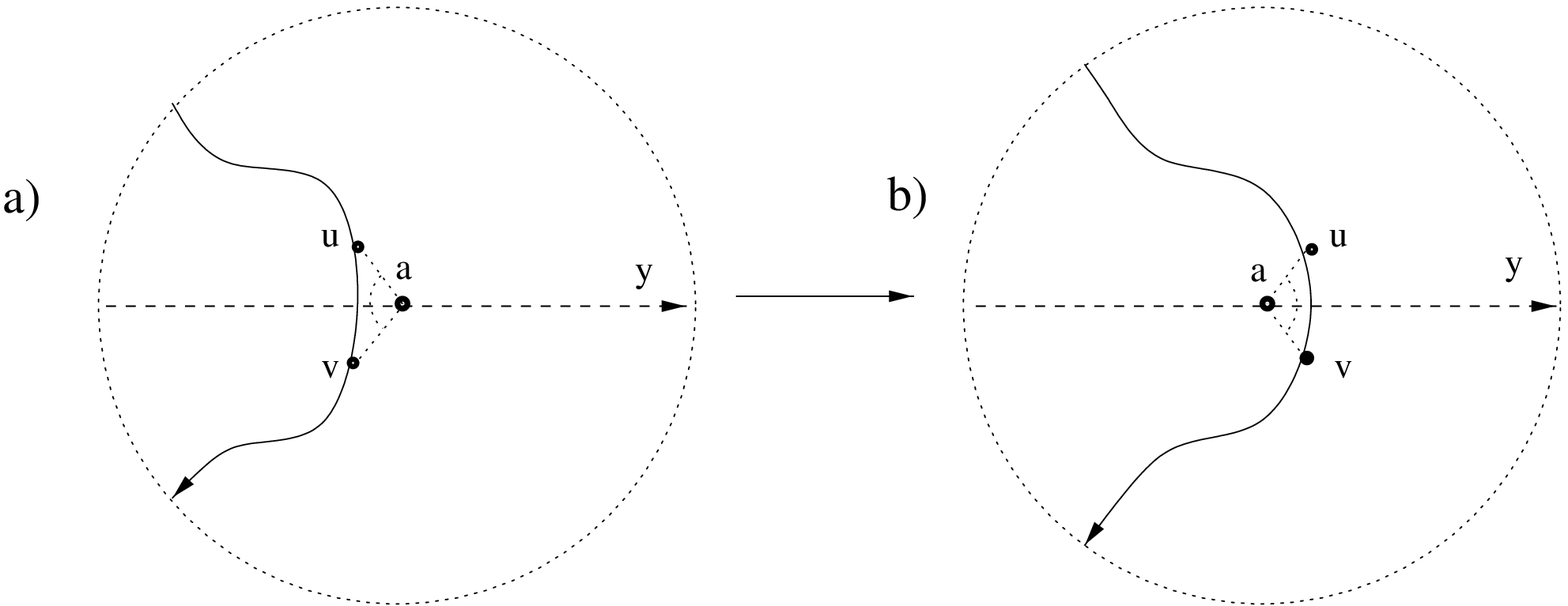}
 \end{center}
\caption{}\label{seif2.fig}
\end{figure}

We can assume that the move proceeds as follows. At the start 
$K$ and $T_a$
intersect along a curve lying in the meridional disk $D$ of $T_a$. The part of
$K$ close to $a$ in $D$ is an arc $C$ of a circle of a very large radius. This
arc is symmetric with respect to the $y$ axis passing through $a$ in $D$.
During the move this arc slides along the $y$ axis through the fiber $a$
(see Figure~\ref{seif2.fig}). 

Clearly two points $u$ and $v$ of $C$ after this move are in the same 
fiber if and only if they are symmetric with respect to the $y$ axis, and the 
angle formed by $v,a,u$ in $D$ is less or equal to $\pi$ and is equal 
to  $\frac{2l\pi}{\mu_a}$ for some $l\in \{1,\dots,\mu_a\}$ 
(see Figure~\ref{seif2.fig}).
They are in the same fiber 
before the move if and only if the angle formed by $u,a,v$ in $D$ 
is less than $\pi$ and is equal to $\frac{2l\pi}{\mu_a}$ for 
some $l\in \{1,\dots,\mu_a\}$ (see Figure~\ref{seif2.fig}).
 
\begin{figure}[htbp]
 \begin{center}
  \epsfxsize 10cm
  \hepsffile{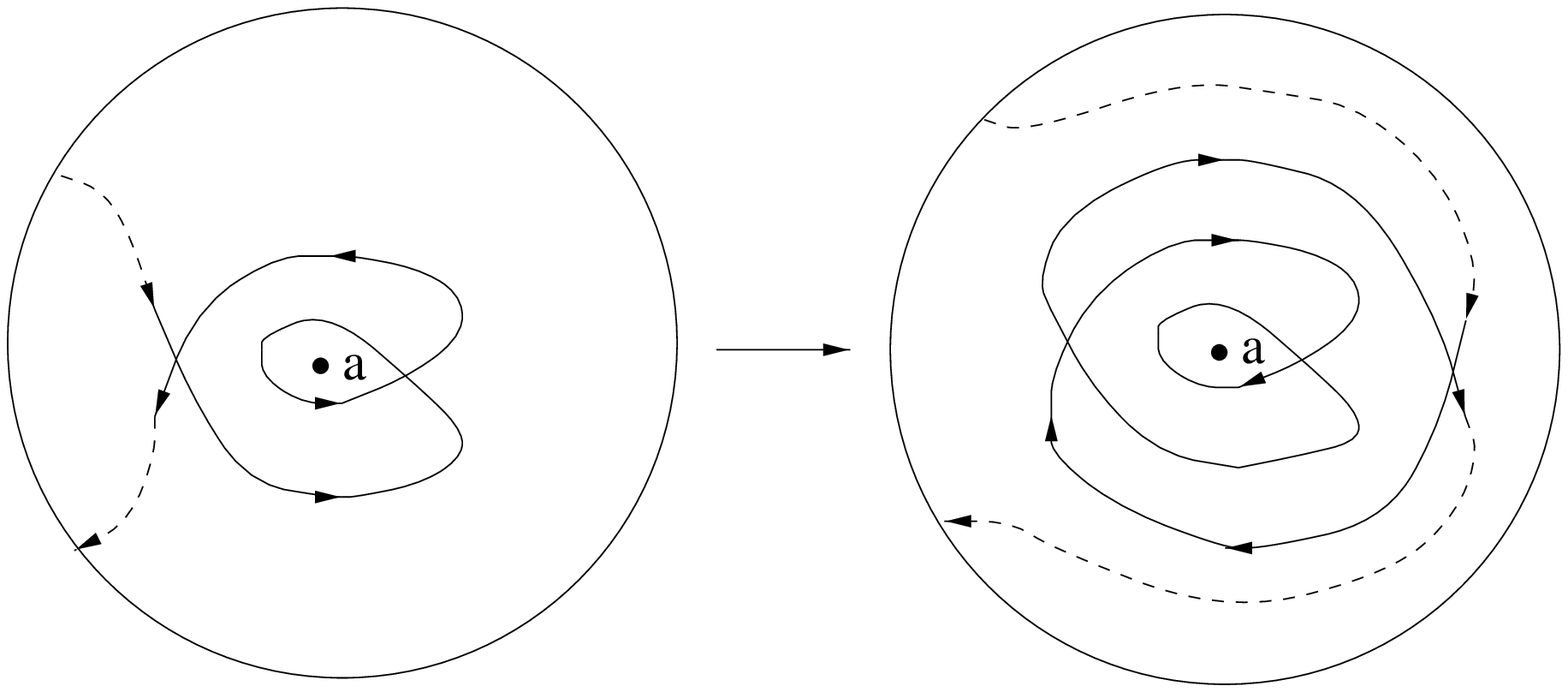}
 \end{center}
\caption{}\label{seif1.fig}
\end{figure}

Consider a double point $v$ of $\pi\big|_D(K)$ that appears 
after the move and
corresponds to the angle $\frac{2l\pi}{\mu_a}$. There is a unique
$k\in N_1(a)$ such that 
$\frac{2\pi\nu_a k}{\mu_a}\mod 2\pi=\frac{2l\pi}{\mu}$. Note that to make
the splitting of $[K]$ into $[\mu^+_v]$ and $[\mu^-_v]$ 
well defined, we do not need the two points of $K$ projecting to $v$ 
to be antipodal in $\pi^{-1}(v)$. This allows one to compare 
these homology classes with $f_a$. For the orientation of $C$ shown in 
Figure~\ref{seif2.fig} one verifies that connecting $v$ to $u$ along the orientation of
the fiber we are adding $k$ fibers $f_a$. (Note that the
factorization we used to define the exceptionally fibered torus was
$\bigl(\bigl(r,\theta\bigr),1\bigr)=
\bigl(\bigl(r,\theta+\frac{2\pi\nu}{\mu} \bigr),0\bigr)$.) 
Thus $[\mu^-_v]=f_a^k$ (see
Figure~\ref{seif1.fig}). From~\ref{sum-property} we know that
$[\mu^+_v][\mu^-_v]=[K]f$. Hence $[\mu^+_v]=[K]f_a^{\mu_a-k}$.

As above, to each double point $v$ of $\pi\big|_D(K)$ before this move
there corresponds $k\in N_2(a)$. For this double point $[\mu^+_v]=f_a^{\mu_a-k}$ and
$[\mu^-_v]=[K]f_a^k$. 

Summing up over the corresponding values of $k$ we see that 
$S_K$ changes by $R^1_a$ under this move. Recall that $R^1_a=0\in H$. Thus 
$S_K$ is invariant under the move.

For the other choice of the orientation of $C$ the value of $S_K$ 
changes by $R^2_a=0\in H$.

Thus $S_K$ is invariant under all elementary isotopies, and this proves
Theorem~\ref{correctSSK}.
\qed

\subsection{Proof of Theorem~\ref{front-shadow}.}\label{pffront-shadow}
\begin{figure}[htbp]
 \begin{center}
  \epsfxsize 10cm
  \hepsffile{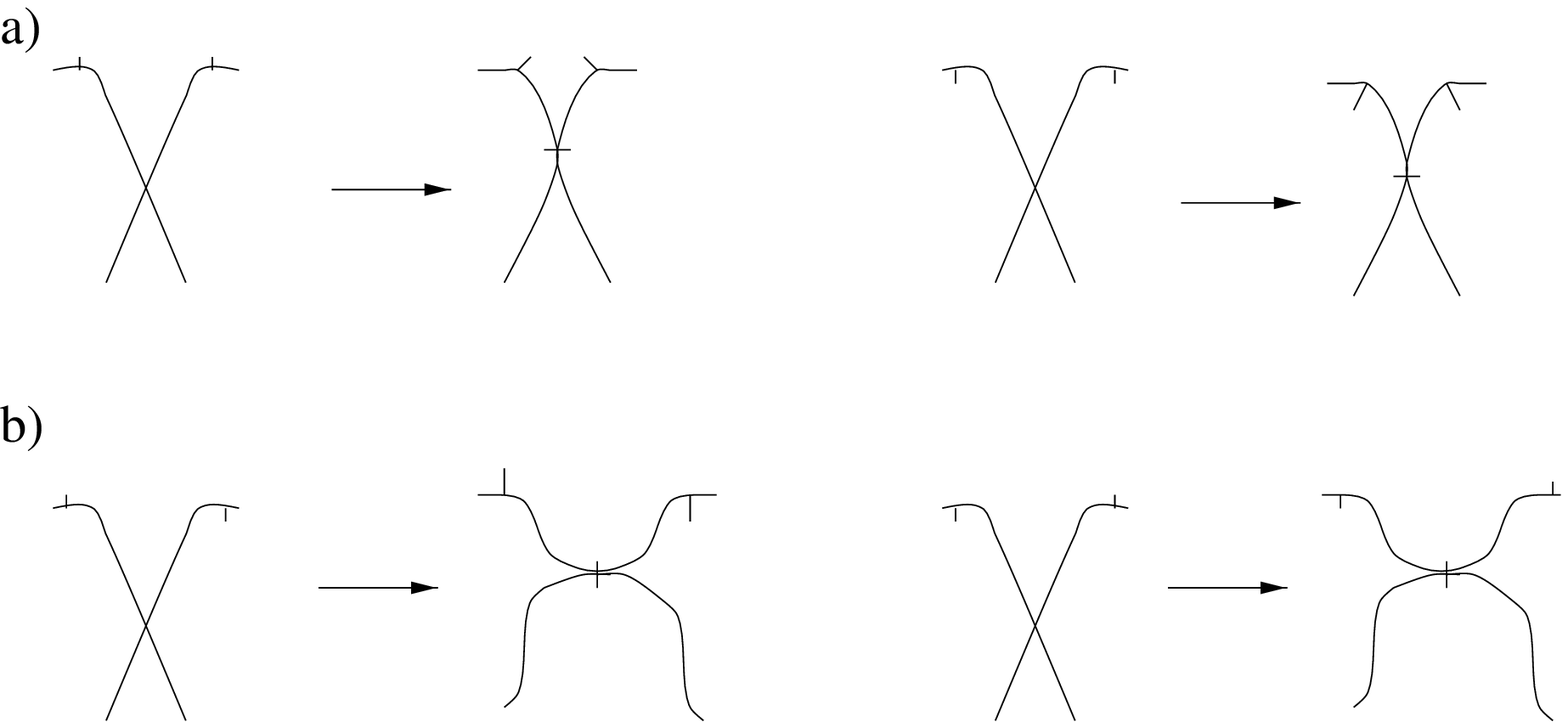}
 \end{center}
\caption{}\label{pic1.fig}
\end{figure}

\begin{figure}[htbp]
 \begin{center}
  \epsfxsize 12cm
  \hepsffile{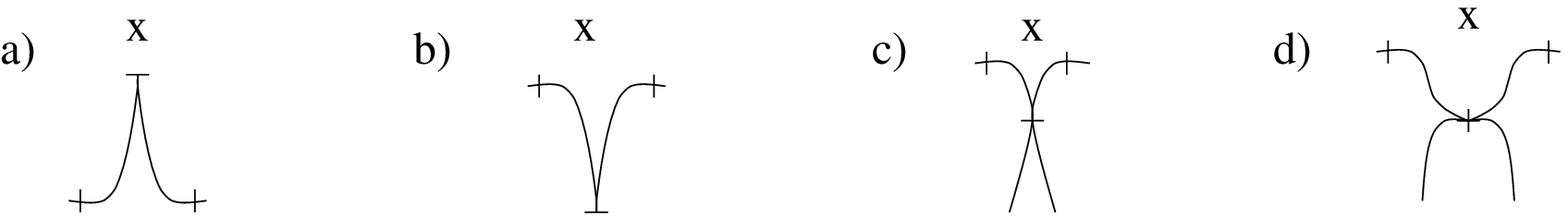}
 \end{center}
\caption{}\label{pic3.fig}
\end{figure}

Deform $L$ in the neighborhoods of all double points of $L$ 
(see Figure~\ref{pic1.fig}), 
so that the two points of the Legendrian knot corresponding 
to the double point of $L$ are antipodal in the fiber. 
After we make the quotient of the fibration by the $\Z_2$-action, 
the projection of the 
deformed $\lambda$ is not a cooriented front anymore but a front 
equipped with a normal field of lines. 
(This corresponds to the factorization $S^1\rightarrow\R P^1$.) 
Using Figure~\ref{pic3.fig} one calculates 
the contributions of different cusps and double points to the total 
rotation number of the line field under traversing the boundary in the 
counter clockwise direction.

These contributions are as follows:
\begin{equation}
\begin{cases}
   1 & \text{for every cusp point pointing inside $X$};\\
  -1 & \text{for every cusp point pointing outside $X$};\\
  -1 & \text{for every double point of the type shown in
Figure~\ref{pic2.fig}c};\\
   0 & \text{for the other types of double points.}
\end{cases}
\end{equation}

To get the contributions to gleams, we divide these numbers by $2$ 
(as we do in the construction of shadows, see~\ref{prelim}).

If the region does not have cusps and double points in its boundary, 
then the obstruction to an extension of the section over $\p X$ 
to $X$ is equal to $\chi(\Int X)$.

This completes the proof of Theorem~\ref{front-shadow}.
\qed

\subsection{Proof of Theorem~\ref{equivalence}.}\label{pfequivalence}
A straightforward verification shows that 
\begin{equation}\label{useful1}
\ind \tilde L_u^+-\ind \tilde L_u^-=
\ind L^+_u-\ind L^-_u-\epsilon(u),
\end{equation} 
\begin{equation}\label{useful2}
\ind \tilde L^+_u+ \ind \tilde L^-_u=\ind L+1,
\end{equation}
and
\begin{equation}\label{useful3}
\ind L^+_u + \ind L^-_u = \ind L
\end{equation}
for any double point $u$ of $L$.

Let us prove~\eqref{rel1}. We write down the formal sums used to
define $S(\lambda)$ and $l_q(\lambda)$ and start to reduce them in a
parallel way as described below. 

We say that a double point is essential if $[\tilde\lambda^+_u]\neq 
[\tilde\lambda^-_u]$. 

For non-essential $u$ we see that the term 
$\bigl([\tilde\lambda^+_u]-[\tilde\lambda^-_u]\bigr)$ in
$S(\lambda)$ is zero. Using~\eqref{useful1} we get that the term
$[\ind L^+_v- \ind L^-_v-\epsilon (v)]_q$ in   
$l_q(\lambda)$ is also zero. 

The index of a wave front coincides with the homology class of its lifting
under the natural identification of $H_1(ST^*F)$ with $\Z$. This fact
and~\eqref{useful2} imply that if we have 
$[\tilde \lambda^+_u]=[\tilde \lambda^-_v]$ for two double points $u$ and $v$, 
then $[\tilde \lambda^-_u]=[\tilde \lambda^+_v]$. 
Hence $\bigl([\tilde \lambda^+_u]-
[\tilde \lambda^-_u]\bigr)=-\bigl([\tilde \lambda^+_v]-
[\tilde \lambda^-_v]\bigr)$,
and these two terms cancel out. Identity~\eqref{useful1} implies that the 
terms $[\ind L^+_u-\ind L^-_u-\epsilon(u)]_q$ and 
$[\ind L^+_v-\ind L^-_v-\epsilon(v)]_q$ also cancel out. 

For similar reasons, if for a double point $u$ the term 
$\bigl ([\tilde 
\lambda^+_u]-[\tilde \lambda^-_u]\bigr)$ is equal to 
$\bigl([\lambda]-f\bigr)$, so that we can simplify $S(\lambda)$ by crossing
out the term and decreasing the coefficient $C^+$ by one. 
Then 
$[\ind L^+_u-\ind L^-_u-\epsilon(u)]_q=[\ind L-1]_q$, and we can simplify 
$l_q(\lambda)$ by crossing out the term and decreasing 
the coefficient $C^+$ by one.

Similarly if the input of  double point $u$ into $S(\lambda)$ is 
$\bigl(e-[\lambda]f\bigr)$, then we reduce the two sums in the parallel way
by crossing out the corresponding terms and decreasing by one the coefficients
$C^-$. 

We make the cancellations described above in both $S(\lambda)$ and
$(l_q(\lambda)-[h]_qh)$ in a parallel way until we can 
not reduce $S(\lambda)$ any more. 
In this reduced form the terms of the form $k_lf^l$ with
$k_l>0$ correspond to 
the terms of type $[\tilde \lambda^+_u]$ for some double points $u$ of $L$.
(The case where a term of this type correspond to cusps is 
treated separately below.) 
Identities~\eqref{useful1} 
and~\eqref{useful2} imply that the contribution of the corresponding 
double points into $l_q(\lambda)$ is $k_l[2l-h-1]_q$. 

In the case where $k_lf^l$ term comes from the cusps and not from the double
points of $L$, one can easily verify that the corresponding input of cusps 
into $(l_q(\lambda)-[h]_qh)$ can still be written as $k_l[2l-h-1]q$.

Thus $l_q(\lambda)=[h]_qh+\sum_{k_l>0} k_l[2l-h-1]_q$, 
and we have proved~\eqref{rel1}.

Let us prove~\eqref{rel2}. 
As above we reduce $S(\lambda)$ and
$(l_q(\lambda)-[h]_qh)$ in a parallel way. Note that the
coefficient at
each $[m]_q$ was positive from the very beginning by the definition
of $l_q(\lambda)$, and it stays positive under the cancellations described  
above. 
After this reduction each term $n_m[m]_q$ 
is a contribution of $n_m$ double points. (The case where it is a
contribution of 
cusps is treated separately as in the proof of~\eqref{rel2}.) 
Let $u$ be one of these double points.
Then from~\eqref{useful1} and~\eqref{useful2} we get
the following system of two equations in variables $\ind \tilde L^+_u$ and 
$\ind \tilde L^-_u$:
\begin{equation}
\begin{cases}
\ind \tilde L^+_u-\ind \tilde L^-_u=m &\\
\ind \tilde L^+_u+\ind \tilde L^-_u=\ind L+1.&
\end{cases}
\end{equation}

Solving the system we get that $[\tilde \lambda
^+_u]=f^{\frac{m+h+1}{2}}$ and $[\tilde \lambda
^-_u]=f^{\frac{h+1-m}{2}}$.

This proves identity~\eqref{rel2} and Theorem~\ref{equivalence}.
\qed

\subsection{Proof of Theorem~\ref{orbifold1}.}\label{pforbifold1}
There are five elementary isotopies of a generic front $L$ on an orbifold
$F$. Four of them are: the birth of two cusps, passing through a
non-dangerous self-tangency point, passing through a triple point,
and passing of a branch through a cusp point. For all possible oriented
versions of these moves a straightforward calculation shows that 
$S(\lambda)\in \frac{1}{2}\Z[H_1(N)]$ is preserved.

\begin{figure}[htbp]
 \begin{center}
  \epsfxsize 11cm
  \hepsffile{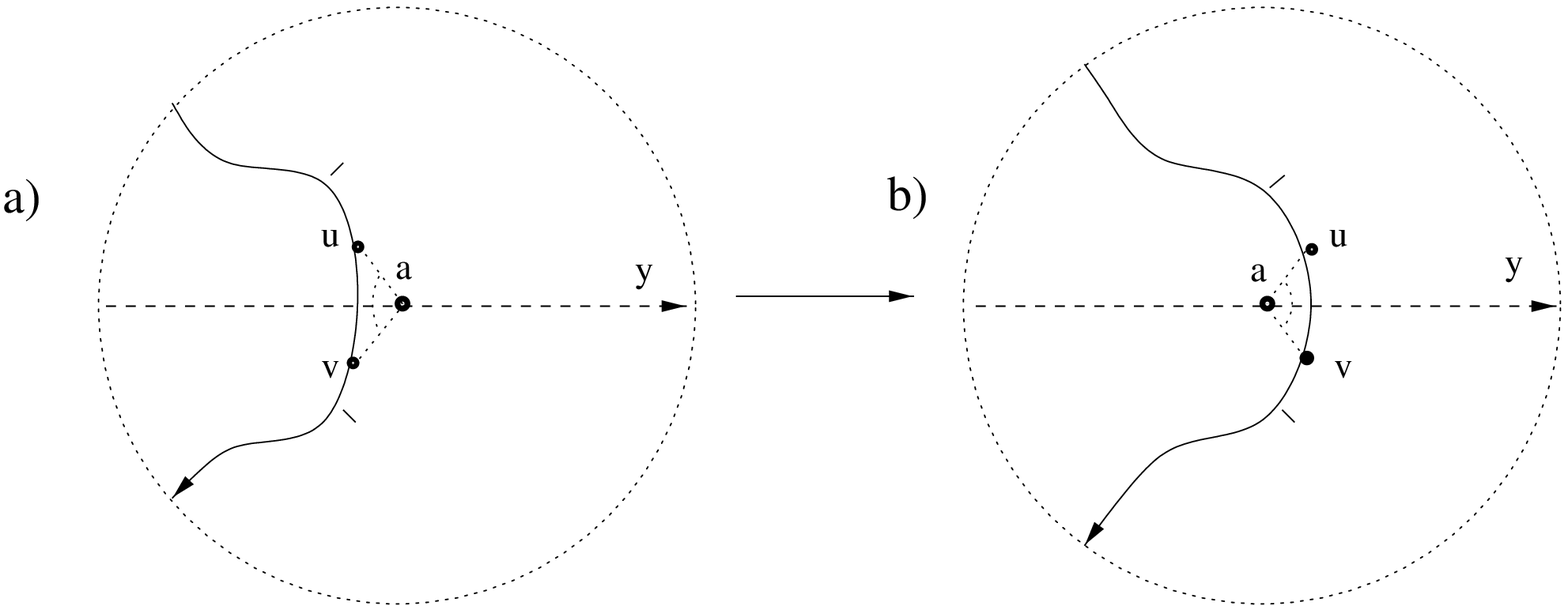}
 \end{center}
\caption{}\label{orbi2.fig}
\end{figure}

The fifth move is more complicated. It corresponds to a generic passing of 
a wave front lifted to $\R^2$ 
through the preimage of a cone point $a$. 
We can assume that this
move is a symmetrization by $G_a$ of the following move. 
The lifted front in the neighborhood of $a$ is an arc $C$ of a circle
of large radius with center at the $y$ axis, and during this move 
this arc slides through $a$ along $y$  (see Figure~\ref{orbi2.fig}).

Clearly after this move points $u$ and $v$ on the arc $C$ turn out to be 
in the same fiber if and only if they are symmetric with respect 
to the $y$ axis, and the angle formed by $v,a,u$ is less or equal to $\pi$ 
and is equal to  $\frac{2k\pi}{\mu_a}$ for some $k\in \{1,\dots,\mu_a\}$ 
(see Figure~\ref{orbi2.fig}). 
We denote the set of such numbers $k$ by 
$\bar N_1(a)=\bigl\{k\in\{1,\dots,\mu_a\}\big|\frac{2k\pi}{\mu_a}
\in(0,\pi]\bigr\}$. 

Two points $u$ and $v$ on the arc $C$ are in the same fiber 
before the move if and only if they are symmetric with respect 
to the $y$ axis, and the angle formed by $u,a,v$ is less than
$\pi$ and is equal $\frac{2k\pi}{\mu_a}$ for 
some $k\in \{1,\dots,\mu_a\}$ (see Figure~\ref{orbi2.fig}).
We denote the set of such numbers $k$ by 
$\bar N_2(a)=\bigl\{k\in\{1,\dots,\mu_a\}\big|\frac{2k\pi}{\mu_a}
\in(0,\pi)\bigr\}$.

\begin{figure}[htbp]
 \begin{center}
  \epsfxsize 10cm
  \hepsffile{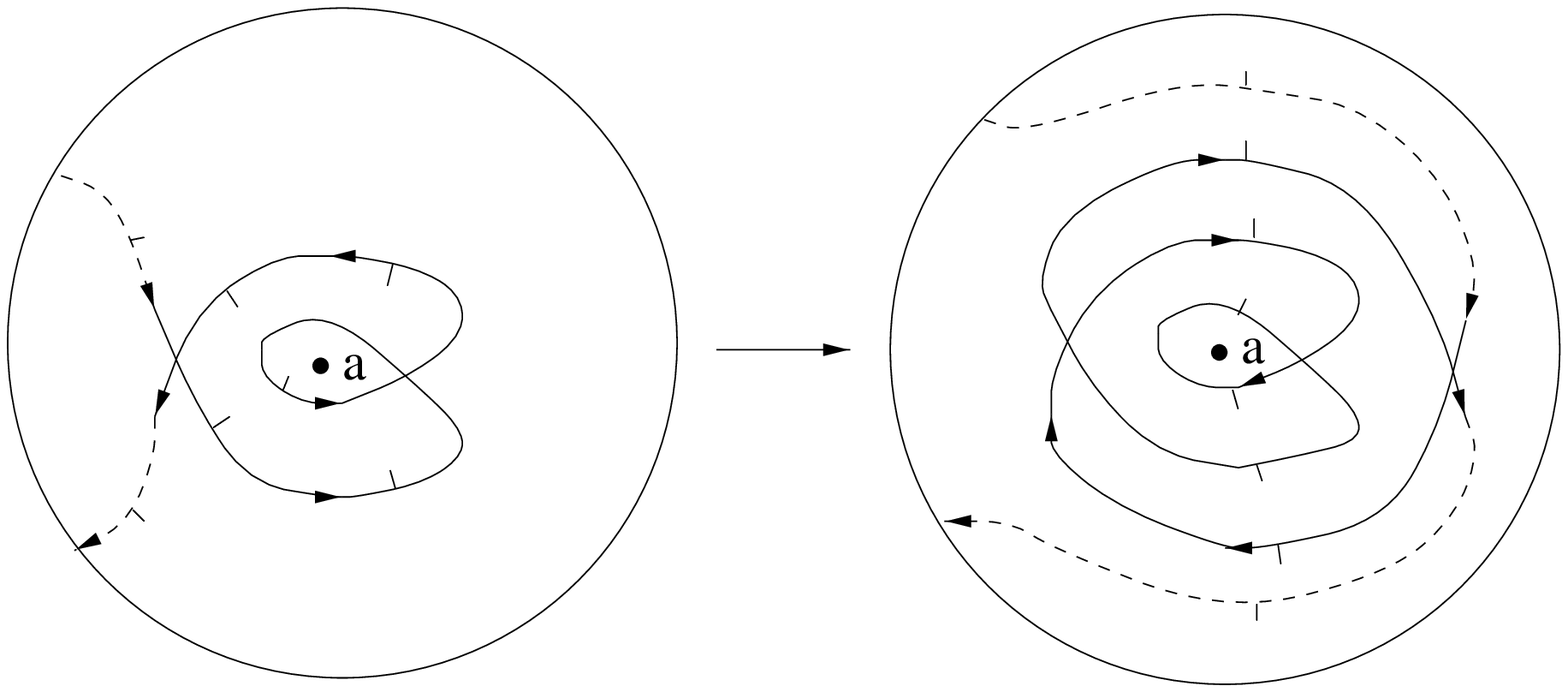}
 \end{center}
\caption{}\label{orbi1.fig}
\end{figure}

The projection of this move for the orientation of $L'$ drawn in
Figure~\ref{orbi2.fig} is shown in Figure~\ref{orbi1.fig}. 

Split the wave front in
Figure~\ref{orbi1.fig} at the double point $v$ (appearing after the move) 
that corresponds to some
$k\in \bar N_1(a)$. 
Then $\tilde\lambda^-_v$ is a front with two positive cusps that
rotates $k$ times around $a$ in the clockwise direction. Hence
$[\tilde\lambda^-_v]=ff_a^{-k}=f_a^{\mu_a-k}$. We know that 
$[\tilde\lambda^+_v][\tilde\lambda^-_v]=[\lambda]f$ and that
$f_a^{\mu_a}=f$. 
Thus
$[\tilde\lambda^+_v]=[\tilde\lambda]f_a^k$.
 
In the same way we verify that if we split the front 
at the double point $v$ 
(existing before the move) that corresponds  to some $k\in \bar N_2(a)$,
then $[\tilde \lambda^+_v]=f_a^k$ and $[\tilde \lambda^-_v]=[K]f_a^{|\mu_a|-k}$. 

Now making sums over all corresponding numbers $k\in \{1,\dots,\mu_a\}$
we get
that under this move $S(\lambda)$ changes by 
\begin{equation}
\bar R^1_a=\sum_{k\in\bar N_1(a)}\bigl([\lambda]f_a^k-f_a^{\mu_a-k}\bigr)-
\sum_{k\in\bar N_2(a)}\bigl(f_a^k-[\lambda]f_a^{\mu_a-k}\bigr).
\end{equation}
A straightforward verification shows that $R^1_a=\bar R^1_a$. (Note that the 
sets $N_1(a)$ and $N_2(a)$ are different from $\bar N_1(a)$ and $\bar N_2(a)$.)

Recall that $R^1_a=0\in J$. Thus $S(\lambda)$ is invariant under the move. 

For the other choice of the orientation of $C$  the value of 
$S(\lambda)$ changes by $R^2_a=0\in J$. 

Hence $S(\lambda)$ is invariant under all elementary
isotopies, and we have proved Theorem~\ref{orbifold1}.
\qed

\end{document}